\documentclass[12pt]{amsart}
\usepackage{amssymb}
\usepackage{amsthm}
\usepackage{amscd}
\usepackage{amsmath}
\usepackage[dvipdfm]{pict2e}
\theoremstyle{plain}
\newtheorem{lemma}{Lemma}[section]
\newtheorem{prop}[lemma]{Proposition}
\newtheorem{thm}[lemma]{Theorem}
\newtheorem{cor}[lemma]{Corollary}
\newtheorem{aplemma}{Lemma~A.\hspace{-1.5mm}}
\newtheorem{approp}{Proposition~A.\hspace{-1.5mm}}
\newtheorem{apthm}{Theorem~A.\hspace{-1.5mm}}
\newtheorem{apcor}{Corollary~A.\hspace{-1.5mm}}
\newtheorem{intthm}{Theorem}

\newtheorem*{problem}{Problem}

\theoremstyle{definition}

\newtheorem{rem}{Remark}
\newtheorem{rema}{Remark}

\newtheorem{remb}{Remark}

\newtheorem{defi}[lemma]{Definition}
\newtheorem{exa}[lemma]{Example}
\newtheorem{aprem}{Remark~A.\hspace{-1.5mm}}
\newtheorem{apdefi}{Definition~A.\hspace{-1.5mm}}
\newcommand{\bde}{\begin{defi}}
\newcommand{\ede}{\end{defi}\vspace{1mm}}
\newcommand{\ble}{\begin{lemma}}
\newcommand{\ele}{\end{lemma}}
\newcommand{\bpr}{\begin{prop}}
\newcommand{\epr}{\end{prop}}
\newcommand{\bt}{\begin{thm}}
\newcommand{\et}{\end{thm}}
\newcommand{\bco}{\begin{cor}}
\newcommand{\eco}{\end{cor}}
\newcommand{\bre}{\begin{rem}}
\newcommand{\ere}{\end{rem}}
\newcommand{\brea}{\begin{rema}}
\newcommand{\erea}{\end{rema}\vspace{1mm}}
\newcommand{\breb}{\begin{remb}}
\newcommand{\ereb}{\end{remb}\vspace{1mm}}
\newcommand{\bex}{\begin{exa}}
\newcommand{\eex}{\end{exa}}
\newcommand{\bpf}{\begin{proof}}
\newcommand{\epf}{\end{proof}\vspace{1mm}}

\newcommand{\bade}{\begin{apdefi}}
\newcommand{\eade}{\end{apdefi}}
\newcommand{\bale}{\begin{aplemma}}
\newcommand{\eale}{\end{aplemma}}
\newcommand{\bapr}{\begin{approp}}
\newcommand{\eapr}{\end{approp}}
\newcommand{\bat}{\begin{apthm}}
\newcommand{\eat}{\end{apthm}}
\newcommand{\baco}{\begin{apcor}}
\newcommand{\eaco}{\end{apcor}}
\newcommand{\bare}{\begin{aprem}}
\newcommand{\eare}{\end{aprem}}

\newcommand{\be}{\begin{enumerate}}
\newcommand{\ee}{\end{enumerate}}
\newcommand{\bcd}{\[\begin{CD}}
\newcommand{\ecd}{\end{CD}\]}
\newcommand{\bit}{\begin{itemize}}
\newcommand{\eit}{\end{itemize}}
\newcommand{\bq}{\begin{quote}}
\newcommand{\eq}{\end{quote}}
\newcommand{\ba}{\begin{array}}
\newcommand{\ea}{\end{array}}
\newcommand{\mcL}{\mathcal{L}}
\newcommand{\mcM}{\mathcal{M}}
\newcommand{\mcN}{\mathcal{N}}
\newcommand{\mcO}{\mathcal{O}}

\newcommand{\mcV}{\mathcal{V}}

\newcommand{\mbD}{\mathbb{D}}

\newcommand{\mbM}{\mathbb{M}}
\newcommand{\mbN}{\mathbb{N}}

\newcommand{\mbP}{\mathbb{P}}
\newcommand{\mbQ}{\mathbb{Q}}

\newcommand{\mbZ}{\mathbb{Z}}
\newcommand{\mfh}{\mathfrak{h}}

\newcommand{\spec}{\mathrm{Spec} \, }

\newcommand{\migi}{\rightarrow}
\newcommand{\longmigi}{\longrightarrow}

\newcommand{\hidari}{\leftarrow}
\newcommand{\longhidari}{\longleftarrow}

\newcommand{\isom}{\stackrel{\sim}{\migi}}

\newcommand{\longisom}{\stackrel{\sim}{\longmigi}}

\newcommand{\migiincl}{\hookrightarrow}
\newcommand{\hidariincl}{\hookleftarrow}
\newcommand{\migisurj}{\twoheadrightarrow}

\newcommand{\mr}{\mathrm}
\newcommand{\hidden}[1]{\,}
\newcommand{\arrangeVa}{\vspace{2mm}}

\begin{document}

\title{On the cuspidalization problem \\ for hyperbolic curves \\  over finite fields}
%\runningheadtitle{cuspidalization problem}

\author{Yasuhiro Wakabayashi}
\date{}
\maketitle
%\markboth{\hfill Yasuhiro Wakabayashi \hfill}{\hfill cuspidalization problem for hyperbolic curves \hfill}

%\footnote{Research Institute for Mathematical Sciences, Kyoto University}
\begin{abstract}
In this paper, we study some group-theoretic  constructions associated to arithmetic fundamental groups of hyperbolic curves over finite fields. One of the main results of this paper asserts that any Frobenius-preserving isomorphism between the geometrically pro-$l$ fundamental groups of hyperbolic curves with one given point removed induces an isomorphism between the geometrically pro-$l$ fundamental groups of the hyperbolic curves obtained by removing other points. Finally, we apply this result to obtain results concerning certain cuspidalization problems for fundamental groups of (not necessarily proper) hyperbolic curves over finite fields.  
\end{abstract}
\tableofcontents 
\section*{Introduction}
\markboth{Yasuhiro Wakabayashi}{cuspidalization problem for hyperbolic curves}

 In the present paper, we consider the following problem:
\begin{problem} \leavevmode\\
 \ \ \ Suppose that we are given a hyperbolic curve over a finite field in which $l$ is invertible.
 Then, given the geometrically pro-$l$ fundamental group of the curve obtained by removing a {\bf specific} point from this hyperbolic curve, is it possible to reconstruct the geometrically pro-$l$ fundamental  groups of the curves obtained by removing other points  which {\bf vary} ``continuously" in a suitable sense?
\end{problem}
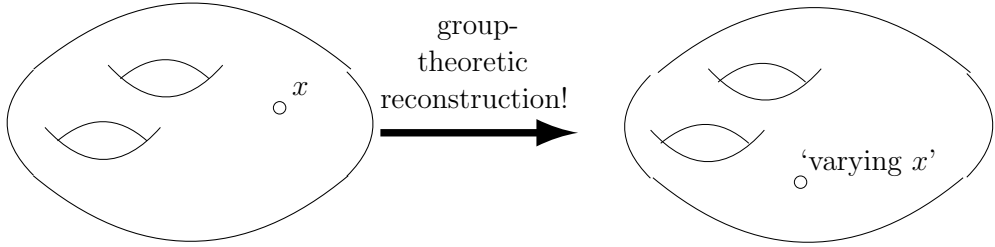
\begin{center}
\begin{picture}(400,150)
\qbezier(20,90)(80,140)(140,90)
\qbezier(20,90)(0,70)(20,50)
\qbezier(20,50)(80,0)(140,50)
\qbezier(140,50)(160,70)(140,90)

\qbezier(380,90)(320,140)(260,90)
\qbezier(380,90)(400,70)(380,50)
\qbezier(380,50)(320,0)(260,50)
\qbezier(260,50)(240,70)(260,90)

\qbezier(25,70)(48,45)(70,70)
\qbezier(30,65)(48,80)(65,65)
\qbezier(50,95)(73,70)(95,95)
\qbezier(55,90)(73,105)(90,90)

\qbezier(265,70)(288,45)(310,70)
\qbezier(270,65)(288,80)(305,65)
\qbezier(290,95)(313,70)(335,95)
\qbezier(295,90)(313,105)(330,90)

\put(120,80){\circle{5}}
\put(125,85){$x$}
\put(330,50){\circle{5}}
\put(330,55){`varying $x$'}
\linethickness{3pt}
\put(160,70){\vector(1,0){80}}

\put(183,110){group-}
\put(175,95){theoretic}
\put(160,80){reconstruction!}

\end{picture}
\end{center}
We shall formulate the above problem mathematically.

Let $l$ be a prime number, $X$  a hyperbolic curve over a finite field $K$ in which $l$ is invertible. For $n$ a positive integer, we denote by $X_n$ the $n$-th configuration space associated to $X$ (hence, $X_1 = X$), and write $\Pi_{X_n}$ for the geometrically pro-$l$ fundamental group of $X_n$.
In the case $n =2$,  the fiber of a projection $X_2 \migi X$ over a $K$-rational point $x \in X$ may be naturally identified with $X \setminus \{ x \}$,  
so we may regard $X_2 \migi X$ as a \textit{continuous family of cuspidalizations of $X$}.
Therefore, the above problem can be formulated as follows (where $Y$  denotes a hyperbolic curve  over a finite field $L$ in which $l$ is also invertible, and we use similar notations for $Y$ to the notations used for $X$):
%-----------------------------------------------------------------------[begin theorem]-------------------
\vspace{3mm}
\begin{intthm} \leavevmode\\
 \ \ \ Let 
\[ \alpha : \Pi_{X \setminus \{ x \} } \stackrel{\sim}{\longmigi} \Pi_{Y \setminus \{ y \} } \]
 be a Frobenius-preserving isomorphism [cf. Definition 3.5] which maps a specific decomposition group $D_{x}$ of $x$  onto a specific decomposition group $D_{y}$ of $y$.
Here, we shall denote by $\overline{\alpha} : \Pi_X \isom \Pi_Y$ (resp., $\overline{D}_x$, $\overline{D}_y$) the isomorphism (resp., as the image of $D_x$ in $\Pi_X$, as the image of $D_y$ in $\Pi_Y$)
obtained by passing to the quotients 
$\Pi_{X \setminus \{ x \} } \migisurj \Pi_X$,  $\Pi_{Y \setminus \{ y \} } \migisurj \Pi_Y$. 

Then there exists an isomorphism
\[ \alpha_2:\Pi_{X_2} \stackrel{\sim}{\longmigi} \Pi_{Y_2} \]
which is uniquely determined up to composition with an inner automorphism  (of either the domain or codomain) by the condition that
it is  compatible with the natural switching automorphisms  up to an inner automorphism  (of either  the domain or codomain) and fits into a commutative diagram
\bcd
\Pi_{X_2} @> \alpha_2 >> \Pi_{Y_2}
\\
@V p_1 VV @VV p_1 V 
\\
\Pi_X @> \overline{\alpha} >> \Pi_Y
\ecd
that induces $\alpha$ by restricting $\alpha_2$ to the inverse images (via the vertical arrows) of  $\overline{D}_x$ and $\overline{D}_y$.

In particular, if $x'$ (resp., $y'$) is a $K$-rational point of $X$ (resp., an $L$-rational point of $Y$), and we assume that the decomposition groups of $x'$, $y'$ correspond via $\alpha$, then we have an isomorphism
\[ \alpha' : \Pi_{X \setminus \{ x' \} } \stackrel{\sim}{\longmigi} \Pi_{Y \setminus \{ y' \} } \]  
(which may not be unique) such that $\alpha$ and $\alpha'$ induce the same isomorphism $\Pi_X \isom \Pi_Y$.
\end{intthm}
\vspace{3mm}
%---------------------------------------------------------------------[end theorem]-------------------

Now let us explain the content of each section briefly. In Section 1,  we recall the notion of the (log) configuration space associated to a hyperbolic curve and review group-theoretic properties of the various fundamental groups associated to such spaces. In particular, the splitting determined by the Frobenius action on the pro-$l$ \'{e}tale fundamental group $\Delta_{X_n}$ of $X_n \times_K \overline{K}$ gives rise to an explicit description of the graded Lie algebra obtained by considering the weight filtration on $\Delta_{X_n}$ (cf. Definition 1.6).  This explicit description will play an essential role in the proof of Theorem A.

 In Section 2, we discuss a certain {\it specific choice}  (among composites with inner automorphisms) of the morphism between geometrically pro-$l$ fundamental groups obtained by switching the two ordered marked points parametrized by the second configuration space.  This choice will play a key role in the proof of Theorem A. 

Section 3 is devoted to proving Theorem A.
Roughly speaking, starting from a given geometrically pro-$l$ fundamental group $\Pi_{X \setminus \{ x \} }$, we reconstruct group-theoretically a suitable topological group, i.e., $\Pi_{X_2}^{\mr{Lie}}$ (cf. Definition 3.1), which contains the geometrically pro-$l$ fundamental group of the second configuration space, by using the explicit description of graded Lie algebra studied in Section 1.
Next, we reconstruct the automorphism on $\Pi_{X_2}^{\mr{Lie}}$ induced by the specific choice of the switching morphism studied in Section 2. 
Finally, we verify that $\Pi_{X_2}$ can be generated, as a subgroup of $\Pi_{X_2}^{\mr{Lie}}$,
by the given fundamental group $\Pi_{X \setminus \{ x \} }$ and the image of this fundamental group via the specific choice of the switiching morphism studied in Section 2; this allows us to 
 reconstruct $\Pi_{X_2}$ as a subgroup of $\Pi_{X_2}^{\mr{Lie}}$.

In Section 4, 
as an application of (a slightly generalized version of) Theorem A, we give a
group-theoretic construction of the cuspidalization of an affine hyperbolic curve $X$ over a finite field at a point {\it ``infinitesimally close"} to the cusp $x$.
That is to say, we give a construction, starting from the geometrically pro-$l$ fundamental group $\Pi_X$ of $X$, of the geometrically pro-$l$ fundamental group $\Pi_{\overline{X}_x^{\mr{log}}}$ of the log scheme obtained by gluing $X$ to a tripod (i.e., the projective line minus three points) at a cusp $x$ of $X$: 
%------------------------------------------------------------------------[begin theorem]------------------
\vspace{3mm}
\begin{intthm} \leavevmode\\
 \ \ \ Let $X$ (resp., $Y$) be an affine hyperbolic curve over a finite field $K$ (resp., $L$),
$x$ a $K$-rational point of $\overline{X} \setminus X$ (\text{resp.}, $y$ an $L$-rational point of  $\overline{Y} \setminus Y$).
Let 
\[ \alpha : \Pi_X \stackrel{\sim}{\longmigi} \Pi_Y \]
 be a Frobenius-preserving isomorphism such that the decomposition groups of $x$ and $y$ (which are well-defined up to conjugacy) correspond via $\alpha$.
Then there exists an isomorphism 
\[ \alpha_{x,y} : \Pi_{\overline{X}_x^{\mr{log}}} \stackrel{\sim}{\longmigi} \Pi_{\overline{Y}_y^{\mr{log}}} \]
which is uniquely determined up to composition with an inner automorphism 
(of either the domain or codomain)
by the condition 
it maps the conjugacy class of the decomposition group of $\tilde{x}$
 to the conjugacy class of decomposition group of $\tilde{y}$, and induce $\alpha$ upon passing to the quotients
$\Pi_{\overline{X}_x^{\mr{log}}} \migisurj \Pi_X$,
$\Pi_{\overline{Y}_y^{\mr{log}}} \migisurj \Pi_Y$.
\end{intthm}
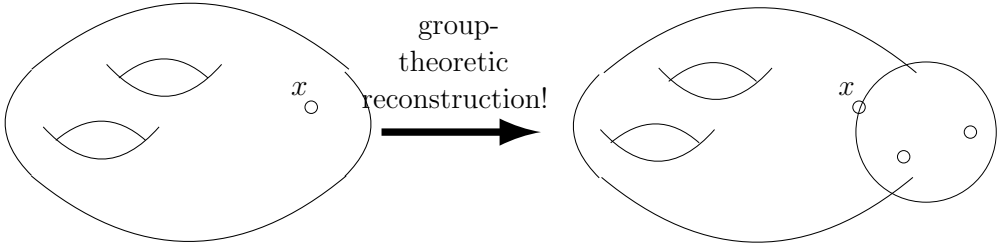
\vspace{3mm}
%------------------------------------------------------------------------[begin theorem]-----------------

\begin{center}
\begin{picture}(400,150)
\qbezier(20,90)(80,140)(140,90)
\qbezier(20,90)(0,70)(20,50)
\qbezier(20,50)(80,0)(140,50)
\qbezier(140,50)(160,70)(140,90)

\qbezier(360,90)(300,140)(240,90)

\qbezier(360,50)(300,0)(240,50)
\qbezier(240,50)(220,70)(240,90)

\qbezier(25,70)(48,45)(70,70)
\qbezier(30,65)(48,80)(65,65)
\qbezier(50,95)(73,70)(95,95)
\qbezier(55,90)(73,105)(90,90)

\qbezier(245,70)(268,45)(290,70)
\qbezier(250,65)(268,80)(285,65)
\qbezier(270,95)(293,70)(315,95)
\qbezier(275,90)(293,105)(310,90)

\put(133,80){\circle{5}}
\put(125,85){$x$}

\put(353,80){\circle{5}}
\put(371,60){\circle{5}}
\put(398,70){\circle{5}}

\put(380,70){\circle{57}}
\linethickness{3pt}
\put(158,70){\vector(1,0){65}}

\put(173,110){group-}
\put(165,95){theoretic}
\put(150,80){reconstruction!}
\put(345,85){$x$}
\end{picture}
\end{center}

At the end of this paper, we consider the cuspidalization problem for (geometrically pro-$l$) fundamental groups of configuration spaces of (not necessarily proper) hyperbolic curves over finite fields (cf. Theorem 4.4):
%---------------------------------------------------------------------[begin theorem]----------------------
\vspace{3mm}
\begin{intthm}
\leavevmode\\
 \ \ \ Let $X$ (resp., $Y$) be a hyperbolic curve over a finite field $K$ (resp., $L$).
Let
\[ \alpha_1 : \Pi_X \stackrel{\sim}{\longmigi} \Pi_Y \]
be a Frobenius-preserving isomorphism. Then  for any $n \in \mbZ_{\geq 0}$, there exists an isomorphism
\[ \alpha_n : \Pi_{X_n} \stackrel{\sim}{\longmigi} \Pi_{Y_n} \]
which is uniquely determined up to composition with an inner automorphism 
(of either the domain or codomain)
by the condition that it is
compatible with the natural respective outer actions of the symmetric group on $n$ letters and makes the diagram
\bcd
\Pi_{X_{n+1}} @> \alpha_{n+1} >> \Pi_{Y_{n+1}} 
\\
@V p_i VV @VV p_i V
\\
\Pi_{X_n} @> \alpha_n >> \Pi_{Y_n}
\ecd
($i = 1, \cdots , n+1$) commute.
\end{intthm}
\vspace{3mm}
%------------------------------------------------------------------------------[end theorem]-------------

Finally, we make a remark on the results in the present paper.
When the curves involved are of genus $\geq 2$, Theorem A may be obtained as an immediate consequence of ~\cite{MZKcusp}, Theorem 3.1; ~\cite{HSHcusp}, Theorem 4.1; ~\cite{HSHcusp}, Corollary 4.1 (i). 
Also, Theorem C is already proved in  ~\cite{MZKcusp} for the case where $n=2$ and $X$ is proper, and in ~\cite{HSHcusp} for the case where $n \geq 3$ and $X$ is proper. 
On the other hand, 
the proof of Theorem A given in the present paper is considerably simpler and more direct than the proofs of ~\cite{MZKcusp} and ~\cite{HSHcusp}.
Indeed, in the present paper, we shall apply Theorem A to give (cf. Theorem C) a substantially simpler proof of ~\cite{HSHcusp}, Theorem 4.1, than the proof given in ~\cite{HSHcusp},
which, moreover, includes, for the first time, the {\it affine case}.

\vspace{5mm}
\hspace{-4mm}{\bf Acknowledgement} \leavevmode\\
 \ \ \ The author would like to express his sincere gratitude to Professors Shinichi Mochizuki and Yuichiro Hoshi for their warm encouragements, suggestions, and many helpful advices, as well
 as to Professor Akio Tamagawa for constructive comments concerning this paper.
The author would like to thank the referee for reading carefully his manuscript and giving him some comments and suggestions.

\vspace{10mm}
%%%%%%%%%%%%%%%%%%%%%%%%%%%%%%%%%%%%%%%%%%%%%%%%%
%-----------------------------------------------------------------------------[begin section]----------
\section*{Notations and Conventions}
%-------------------------------------------------------------------------------------------------------
\vspace{6mm}
\hspace{-5mm} {\bf Numbers:} \leavevmode\\
 \ \ \ We shall denote by $\mbQ$ the field of
 \textit{rational numbers}, by $\mathbb{Z} $ the ring of \textit{rational integers},
and by $\mathbb{N} \subseteq \mbZ \,\bigl( \mr{resp.},\,\mbZ_{ \geq a} \subseteq \mbZ \bigr)$ the additive submonoid of integers $\textit{n} \ge 0$ $\bigl($resp., the subset of integers $n \geq a$ for $a \in \mbZ \bigr)$.
If \textit{l} is a prime number, then $\mathbb{Z}_l \ \bigl( $resp.,$ \ \mathbb{Q}_l \bigr)$ denotes the \textit{l-adic completion} of $\mathbb{Z} \ \bigl( $resp.,$ \ \mathbb{Q} \bigr)$.

%-----------------------------------------------------------------------------------------------
\vspace{4mm}
\hspace{-5mm} 
{\bf Topological Groups:}\leavevmode\\
 \ \ \ For an arbitrary Hausdorff topological group $G$, the notation
\[ G^{\mr{ab}} \]
will be used to denote the \textit{abelianization} of $G$, i.e., the quotient of $G$ by the closed subgroup of \textit{G} topologically generated by the commutators of $G$.
We shall say that $G$ is {\it slim} if each open subgroup of $G$ is center-free.

For each closed subgroups $H$ of $G$, let us write
\[N_G (H):=\bigl\{g \in G \mid g\cdot H\cdot g^{-1}=H \bigr\} \]
for the \textit{normalizer} of $H$ in $G$.
We shall say that a closed subgroup $H \subseteq G$ is \textit{normally terminal} in $G$ if the normalizer
$N_G(H)$ is equal to $H$.

We shall write $\mr{Aut}(G)$ for the group of automorphisms of the topological group $G$,
$\mr{Inn} :G \migi \mr{Aut}(G)$ for the homomorphism  obtained by letting $G$ act on $G$ by inner automorphisms, and $\mr{Out}(G) := \mr{Aut}(G)/\mr{Inn}(G)$.

If $G'$ is a topological group, then one may define an equivalence relation $\sim$ on $\mr{Hom}(G', G)$, the set of continuous group homomorphisms from $G'$ to $G$, by 
\[   f_1 \sim f_2 \ \Longleftrightarrow \ \exists g \in G: f_1 = \mr{Inn}(g) \circ f_2  \]
where $f_1$, $f_2 \in \mr{Hom}(G', G)$.
%, and $\mr{Inn}$ denotes the homomorphism $G \migi \mr{Aut}(G)$ obtained by letting $G$ act on $G$ by inner automorphisms.
We shall refer to an element of the quotient set $\mr{Hom}(G', G)/\sim$ of $\mr{Hom}(G', G)$ by $\sim$ as an {\it outer homomorphism}.
Note that $\sim$ is compatible with composition of homomorphisms, so composition of outer homomorphisms make sense.

If $G$ is a center-free, then the natural sequence
 \[ 1 \longmigi G \stackrel{\mr{Inn}}{\longmigi} \mr{Aut}(G) \longmigi \mr{Out}(G) \longmigi 1 \]
 is evidently exact.  If the profinite group $G$ is topologically finitely generated, then the groups $\mr{Aut}(G)$, $\mr{Out}(G)$ are naturally endowed with a profinite topology, and the above sequence may be regarded as an exact sequence of profinite groups. 

%If $G$ is a center-free, then $\mr{Inn} : G \migi \mr{Aut}$ we have a natural exact sequence
% \[ 1 \longmigi G \stackrel{\mr{Inn}}{\longmigi} \mr{Aut}(G) \longmigi \mr{Out}(G) \longmigi 1, \]
 %where $\mr{Aut}(G)$ denotes the group of automorphisms of the topological group $G$; the injective (since $G$ is center-free) homomorphism $\mr{Inn} :G \migi \mr{Aut}(G)$ is obtained by letting $G$ act on $G$ by inner automorphisms; $\mr{Out}(G)$ is defined so as to render the sequence exact. If the profinite group $G$ is topologically finitely generated, then the groups $\mr{Aut}(G)$, $\mr{Out}(G)$ are naturally endowed with a profinite topology, and the above sequence may be regarded as an exact sequence of profinite groups. 

 If, moreover,  $J \migi \mr{Out}(G)$ is a homomorphism of groups, then we shall write
 \[ G \stackrel{\mr{out}}{\rtimes} J := \mr{Aut}(G) \times_{\mr{Out}(G)} J \]
 for the ``\textit{outer semi-direct product of $J$ with $G$}". Thus, we have a natural exact sequence
 \[ 1 \longmigi G \longmigi G \stackrel{\mr{out}}{\rtimes} J \longmigi J \longmigi 1 . \]  
It is verified (cf. ~\cite{HSHcusp}, Lemma 4.10) that if  an automorphism $\phi$ of $G \stackrel{\mr{out}}{\rtimes} J$ preserves the subgroup $G \subseteq G \stackrel{\mr{out}}{\rtimes} J$ and induces  the identity morphism on $G$ and the quotient $J$, then $\phi$ is the identity morphism of $G \stackrel{\mr{out}}{\rtimes} J$. 
%----------------------------------------------------------------------------------------------
\vspace{4mm}

\hspace{-5mm}{\bf  Log schemes:}\leavevmode\\
 \ \ \  Basic references for the notion of \textit{log scheme} are ~\cite{KATO} and ~\cite{ILL}.
In this paper, log structures are always considered on the \'{e}tale sites of schemes.
 For a log scheme $X^{\mr{log}}$, we shall denote by $X$ (resp.,\,$\mcM_X$) the underlying scheme of $X^{\mr{log}}$ (resp.,\,the sheaf of monoids defining the log structure of $X^{\mr{log}}$). 
 Let $X^{\mr{log}}$ and $Y^{\mr{log}}$ be log schemes, and $f^{\mr{log}} :X^{\mr{log}} \migi Y^{\mr{log}}$ a morphism of log schemes. 
Then we shall refer to the quotient of $\mcM_X$ by the image of the morphism $f^{*} \mcM_Y \migi \mcM_X$ induced by $f^{\mr{log}}$ as the \textit{relative characteristic sheaf} of $f^{\mr{log}}$. 
Moreover, we shall refer to the relative characteristic sheaf of the morphism $X^{\mr{log}} \migi X$ (where, by abuse of notation, we write $X$ for the log scheme obtained by equipping $X$ with the trivial log structure) induced by the natural inclusion
$ \mcO^* \migiincl \mcM_X $ as the \textit{characteristic sheaf} of $X^{\mr{log}}$.

 We shall say that a log scheme $X^{\mr{log}}$ is \textit{fs} if $\mcM_X$ is a sheaf of integral monoids (cf. ~\cite{ILL}, \S\,1.1), and locally for the \'{e}tale topology, has a chart modeled on a finitely generated and saturated monoid (cf. ~\cite{ILL}, \S\,1.1).
If $X^{\mr{log}}$ is \textit{fs}, then, for $n$ a nonnegative integer, we shall refer to as the \textit{$n$-interior} of $X^{\mr{log}}$ the open subset of $X$ on which the associated sheaf of 
group envelopes (cf. ~\cite{ILL}, \S\,1.1)
%groupifications
 of characteristic sheaf of $X^{\mr{log}}$ is of $rank \leq n$.  Thus, the $0$-interior of $X^{\mr{log}} $ is often referred to simply
as the {\it interior} of $X^{\mr{log}}$.

%---------------------------------------------------------------------------------------------------------
\vspace{4mm}
\hspace{-5mm}{\bf Curves:}\leavevmode\\
 \ \ \ Let $ f : X \migi S $ be a morphism of schemes. 
Then we shall say that $f$ is a \textit{family of curves of type (g,r)} 
if it factors $X \migiincl \overline{X} \migi S$ as the composite of an open immersion $X \migiincl \overline{X}$ whose image is the complement $\overline{X} \setminus D$ of a relative divisor $D \subseteq \overline{X}$  which is finite \'{e}tale over $S$ of relative degree $r$, 
and a morphism $\overline{X} \migi S$ which is proper, smooth, and geometrically connected, and whose geometric fibers are one-dimensional of genus $g$.
We shall refer to $\overline{X}$ as the \textit{compactification} of $X$.

We shall say that $f$ is a \textit{family of hyperbolic curves (resp.,\,tripod)} 
if $f$ is a family of curves of type ($g,r$) such that ($g,r$) satisfies $2g-2+r > 0$ (resp.,\,$(g,r)=(0,3)$ and the relative divisor $D$ is split over $S$).

We shall denote by 
\[ \overline{\mcM}_{g,[r]+s} \]
the moduli stack of $r+s$-pointed stable curves of genus $g$ for which $s$ sections are equipped with an ordering.  This moduli stack may be obtained as the quotient of the moduli stack of ordered $(r+s)$-pointed stable curves of genus $g$ (cf. ~\cite{KNUD} for an exposition of the theory of such curves) by a suitable symmetric group action on $r$ letters.  We shall denote by $\overline{\mcM}^{\mr{log}}_{g,[r]+s}$ the log stack obtained by equipping $\overline{\mcM}_{g,[r]+s}$ with the log structure associated to the divisor with normal crossings which parametrizes singular curves.
%----------------------------------------------------------------------------------------------------------
\vspace{4mm}

\hspace{-5mm}{\bf Fundamental Groups:}\leavevmode\\
 \ \ \ A basic reference for the notion of \textit{Kummer \'{e}tale covering} is ~\cite{ILL}.
  For a locally Noetherian, connected scheme $X$ (resp., a locally Noetherian, connected, fs log scheme $X^{\mr{log}}$) equipped with a geometric point
$\overline{x} \migi X$ (resp., log geometric point $\tilde{x}^{\mr{log}} \migi X^{\mr{log}}$),
we shall denote by $\pi_1(X, \overline{x})$ (resp., $\pi_1(X^{\mr{log}}, \tilde{x}^{\mr{log}})$) the \'{e}tale fundamental group of $X$ (resp., logarithmic fundamental group of $X^{\mr{log}}$). Since one knows that the \'{e}tale and logarithmic fundamental groups are determined up to inner automorphisms independently of the choice of basepoint, we shall omit the basepoint, 
and write $\pi_1(X)$ (resp., $\pi_1(X^{\mr{log}})$ ).

 For a scheme $X$ (resp., fs log scheme $X^{\mr{log}}$) which is geometrically connected and of finite type over a field $K$ in which a prime number $l$ is invertible, 
we shall refer to the quotient $\Pi_X$ of $\pi_1(X)$ (resp., the quotient $\Pi_{X^{\mr{log}}}$ of $\pi_1(X^{\mr{log}})$) by the closed normal subgroup obtained as the kernel of the natural projection from $\pi_1(X \times_K \overline{K})$ (resp., $\pi_1(X^{\mr{log}} \times_K \overline{K})$)
(where $ \overline{K}$ is a separable closure of $K$) 
to its maximal pro-$l$ quotient $\Delta_X$ (resp., $\Delta_{X^{\mr{log}}}$) as the \textit{geometrically pro-l \'{e}tale fundamental group} of $X$ (resp., \textit{geometrically pro-$l$ logarithmic fundamental group} of $X^{\mr{log}}$). 
Thus, (if we write $G_K$ for the Galois group of a separable closure of $K$ over $K$, then) we have a natural exct sequence 
\[ 1 \longmigi \Delta_X \longmigi \Pi_X \longmigi G_K \longmigi 1 \]
\[ (\text{resp.}, 1 \longmigi \Delta_{X^{\mr{log}}} \longmigi \Pi_{X^{\mr{log}}} \longmigi G_K \longmigi 1). \]
\leavevmode\\ 
Note that if the log structure of $X^{\mr{log}}$ is trivial, then we have natural isomorphisms $\Delta_{X^{\mr{log}}} \isom \Delta_{X}, \Pi_{X^{\mr{log}}} \isom \Pi_{X}$.

If $K$ is finite, then write $G_K^\dagger \subseteq G_K$ for the (unique) \textit{maximal pro-$l$ subgroup of $G_K$} (so $G_K^\dagger \cong \mbZ_l$).
Also, for a profinite group $\Pi$ over $G_K$, we shall use the notation
\[ \Pi^\dagger := \Pi \times_{G_K} G_K^\dagger \subseteq \Pi . \]
and refer to it as the \textit{restricted pro-$l$ group} of $\Pi$. 
\vspace{10mm}
%%%%%%%%%%%%%%%%%%%%%%%%%%%%%%%%%%%%%%%%%%%%%%%%
%%%%%%%%%%%%%%%%%%%%%%%%%%%%%%%%--[ begin  section1]---%%%%%%
\section{Fundamental groups of (log) configuration spaces }\leavevmode\\
 \ \ \ The purpose of this section is to recall the notion of the (log) configuration space associated to
a curve and review group-theoretic properties of the various fundamental groups associated to such spaces.
\\[3mm]
 \ \ \ Let $l$ be a prime number, $K$ a finite field in which $l$ is invertible, $\overline{K}$ a separable closure of $K$, where we shall denote by $G_K$ the Galois group of $\overline{K}$ over $K$, and $X$ a hyperbolic curve over $K$ of type ($g,r$).
%\\
%-------------------------------------------------------------------------[begin definition]--------------
\begin{defi}\leavevmode\\
\be
\vspace{-7mm}
\arrangeVa\item[(i)]
For  $n \in \mbZ_{\geq 1}$, write $X^{\times n}$ for the fiber product of $n$ copies of $X$ over $K$.
We shall denote by 
\[ X_n \bigl( \subseteq X^{\times n} \bigr) \]
 the \textit{$n$-th configuration space associated to $X$}, i.e., 
the scheme which represents the open subfunctor
\[ S \mapsto \bigl\{(f_1,\dotsm ,f_n) \in X^{\times n}(S) \ \big| \ 
f_i \neq f_j \ \mr{if} \ i \neq j \bigr\} \]
of the functor represented by $X^{\times n}$.
\arrangeVa\item[(ii)]
Let us denote by $\overline{X}^{\mr{log}}_n$ the \textit{$n$-th log configuration space associated to X} (cf. ~\cite{MZKTAM}), i.e.,
\[ \overline{X}^{\mr{log}}_n:=\spec K \times_{{\overline{\mcM}}^{\mr{log}}_{g,[r]}} {\overline{\mcM}}^{\mr{log}}_{g,[r]+n} \]
--- where the (1-)morphism $\spec K \migi {\overline{\mcM}}^{\mr{log}}_{g,[r]}$ is the classifying
morphism determined by the curve $X \migi \spec K$, and the (1-)morphism 
${\overline{\mcM}}^{\mr{log}}_{g,[r]+n} \migi {\overline{\mcM}}^{\mr{log}}_{g,[r]}$
is obtained by forgetting the ordered $n$ marked points of the tautological family of curves over $\overline{\mcM}^{\mr{log}}_{g,[r]+n}$ ---.
In the following, for simplicity, we shall write $\overline{X}^{\mr{log}}$ for $\overline{X}_1^{\mr{log}}$.
\ee
\end{defi}
%------------------------------------------------------------------------------------------------------------
%\vspace{2mm}
\bpr\leavevmode\\
\be
\vspace{-7mm}
\arrangeVa\item[(i)]
The $0$-interior (cf. \S\,0) of the log scheme $\overline{X}^{\mr{log}}_n$ is naturally isomorphic to the $n$-th configuration space $X_n$ associated to X.
\arrangeVa\item[(ii)]
The log scheme $\overline{X}^{\mr{log}}_n$ is log regular and its underlying scheme is connected and regular.
\arrangeVa\item[(iii)]
The projection $p_k^{\mr{log}} : \overline{X}_n^{\mr{log}} \migi \overline{X}_{n-1}^{\mr{log}} $, induced from the (1-)morphism 
${\overline{\mcM}}^{\mr{log}}_{g,[r]+n} \migi {\overline{\mcM}}^{\mr{log}}_{g,[r]+n-1}$
 obtained by forgetting the $k$-th ($k = 1, \cdots , n$) ordered points of the tautological family of curves over $\overline{\mcM}^{\mr{log}}_{g,[r]+n}$, is log smooth\,(cf. \S\,0) and its underlying morphism of schemes is the natural projection $p_k : X_n \migisurj X_{n-1}$ obtained by forgetting the $k$-th factor, and hence, is flat, 
 and has connected and reduced fibers over the geometric points of $X_{n-1}$.
% geometrically connected, and geometrically reduced.
\ee
\epr 
\bpf
See, for example, ~\cite{HSHcusp}, Proposition 2.2.
\epf
%----------------------------------------------------------------------------------------------------------
 \bde \leavevmode\\
 \ \ \  We shall denote (cf. \S\,0) by 
 \[  \Pi_{X_n} \ \  (\text{resp.},  \Delta_{X_n} ) \]
 the geometrically pro-$l$ \'{e}tale fundamental group of $X_n$ (resp., $X_n \times_K \overline{K}$), and 
 \[ \Pi_{\overline{X}_n^{\mr{log}}} \  \ (\text{resp.,} \Pi_{\overline{X}^{\mr{log} \times n }}) \]
 the geometrically pro-$l$ log fundamental group of $\overline{X}_n^{\mr{log}}$ (resp.,  the fiber product $\overline{X}^{\mr{log} \times n}$ of $n$ copies of $\overline{X}^{\mr{log}}$ over $K$).  
Moreover, we shall denote (cf. \S\,0) by
\[  \Pi_{X_n}^\dagger, \ \   \Delta_{X_n}^\dagger (\cong \Delta_{X_n}), \ \  \Pi_{\overline{X}_n^{\mr{log}}}^\dagger, \  \ \Pi_{\overline{X}^{\mr{log} \times n }}^\dagger \]
respective restricted geometrically pro-$l$ groups.

Also we shall  write
\[ p_k^\Delta :\Delta_{X_n} \migisurj \Delta_{X_{n-1}}, \ \  p_k^\Pi :\Pi_{X_n} \migisurj \Pi_{X_{n-1}} \] 
for the morphisms induced by the projection $p_k \times_K \overline{K} : X_n \times_K \overline{K} \migisurj X_{n-1} \times_K \overline{K}$, $p_k : X_n \migisurj X_{n-1}$ obtained by forgetting the $k$-th factor
(these morphisms of profinite groups are only defined up to conjugacy in the absence of appropriate choices of basepoints of respective schemes)
and write

\[ i_k^\Delta : \Delta^k_{X_{n / n-1}}\migiincl \Delta_{X_n} , \ \ i_k^{\Delta' }: \Delta^k_{X_{n / n-1}}\migiincl \Pi_{X_n} \]
for the kernels of the surjections
 $ p_k^\Delta :\Delta_{X_n} \migisurj \Delta_{X_{n-1}}$,
$p_k^\Pi :\Pi_{X_n} \migisurj \Pi_{X_{n-1}}$. 
Then we have exact sequences
\[ 1 \longmigi \Delta_{X_n} \longmigi \Pi_{X_n}^{(-)} \longmigi G_K^{(-)} \longmigi 1 \]
\[ 1 \longmigi \Delta^k_{X_{n/n-1}} \stackrel{i_k^\Delta}{\longmigi} \Delta_{X_n} \stackrel{p_k^\Delta}{\longmigi} \Delta_{X_{n-1}} \longmigi 1 \]
\[ 1 \longmigi \Delta^k_{X_{n/n-1}} \stackrel{i_k^{\Delta'}}{\longmigi} \Pi_{X_n}^{(-)} \stackrel{p_k^{\Pi^{(-)}}}{\longmigi} \Pi_{X_{n-1}}^{(-)} \longmigi 1 \]
 --- where the symbol $(-)$ denotes either the presence or absence of ``$\dagger$'' ---.

Also, we have a square diagram
\bcd
\Pi_{X_{n-1}}^{( - )} @< p_k^{\Pi^{(-)}} << \Pi_{X_n}^{( - )} @>>> \overbrace{\Pi_{X}^{(-)} \times_{G_K^{(-)}} \cdots \times_{G_K^{(\_)}} \Pi_{X}^{(-)}}^{n}
\\
@VVV @VVV @VVV
\\
\Pi_{\overline{X}^{\mr{log}}_{n-1}}^{( -)} @< << \Pi_{\overline{X}^{\mr{log}}_n}^{( -)} @>>> \Pi_{\overline{X}^{\mr{log} \times n}}^{(-)},
\ecd
--- which can be made commutative without conjugate-indeterminacy by choosing compatible base points --- arising from a natural commutative diagram
\bcd
X_{n-1} @< p_k << X_n @>>> X^{n \times}
\\
@VVV @VVV @VVV
\\
\overline{X}_n^{\mr{log}} @< p_k^{\mr{log}} << \overline{X}_n^{\mr{log}} @>>> \overline{X}^{\mr{log} \times n} . 
\ecd
Then, it follows from Proposition 1.2 (i), (ii) together with the log purity theorem (cf. ~\cite{ILL}, ~\cite{MZKext}) that the three vertical homomorphisms are isomorphisms. In the following, we shall identify $\Pi_{X_n}^{(-)}$ with $\Pi_{\overline{X}_n^{\mr{log}}}^{(-)}$, $\Pi_{\overline{X}^{\mr{log} \times n }}^{(-)} $ with $\overbrace{\Pi_{X}^{(-)} \times_{G_K^{(-)}} \cdots \times_{G_K^{(\_)}} \Pi_{X}^{(-)}}^{n}$ and the surjection $p_k^\Pi : \Pi_{X_n} \migi \Pi_{X_{n-1}}$
with the surjection $\Pi_{\overline{X}_n^{\mr{log}}}^{(-)} \migi \Pi_{\overline{X}_{n-1}^{\mr{log}}}^{(-)}$
by means of these specific isomorphisms.
\ede
%----------------------------------------------------------------------------[begin proposition]-----------
\bpr\leavevmode\\
\be
\vspace{-7mm}
\arrangeVa\item[(i)]
 $\Delta^k_{X_{n/n-1}}$ may be naturally identified with the maximal pro-$l$ quotient of the \'{e}tale fundamental group of a geometric fiber of the projection morphism $p_k:X_n \migi X_{n-1} $.
\arrangeVa\item[(ii)]
The images of the $i_k^\Delta : \Delta_{X_{n/n-1}}^k \migi \Delta_{X_n}$, where $k= 1, \cdots, n$, generate $\Delta_{X_n}$.

\arrangeVa\item[(iii)]
The profinite groups $\Delta_{X_n}$, $\Delta^k_{X_{n/n-1}}$, $\Pi_{X_n}^\dagger$, $\Pi_{X^{\times n}}^\dagger$ are slim (cf. \S\,0)
% (i.e., every open subgroup of each profinite group is center-free).
\ee
\epr
%------------------------------------------------------------------------[begin proof]-------------------
\bpf Assertion (i) follows from ~\cite{MZKTAM}, Proposition 2.2, or ~\cite{STX}, Proposition 2.3. 
Assertions (ii) and (iii) follow from induction on $n$, together with the exact sequence
\[ 1 \longmigi \Delta^n_{X_{n/n-1}} \stackrel{i_n^\Delta}{\longmigi} \Delta_{X_n} \stackrel{p_n^\Delta}{\longmigi} \Delta_{X_{n-1}} \longmigi 1 \]
displayed in Definition 1.3.
Indeed, with regard to (ii), $\Delta_{X_{n/n-1}}^k$ maps to $\Delta_{X_{n-1/n-2}}^k$ (for $k = 1, \cdots n-1$) via $p_n^\Delta : \Delta_{X_n} \migi \Delta_{X_{n-1}}$, and it is verified that this map 
$\Delta_{X_{n/n-1}}^k\migi \Delta_{X_{n-1/n-2}}^k$ is surjective by regarding it as the morphism induced by an open immersion between the hyperbolic curves that arise as geometric fibers of the projection morphisms involved.
With regard to (iii), the slimness of $\Delta_X$ is well-known (cf., e.g., ~\cite{MZKhyp}, Lemma 1.3.10); the slimness of $\Pi_X^\dagger$ follows from the fact that the character of $G_K^\dagger$ arising from the determinant of $\Delta_X^{\mr{ab}}$ coincides with some positive power of the cyclotomic character; the other statements follow from the fact that an extension of slim profinite groups is itself slim.
\epf
%---------------------------------------------------------------------------[end proposition]--------------
Next, we recall from ~\cite{MZKcusp}, \S\,3, the theory of the weight filtration of fundamental groups
and the associated graded Lie algebra.
%------------------------------------------------------------------------[begin definition]--------------
\bde
\leavevmode\\
 \ \ \ Let $l$ be a prime number;
$G$, $H$, $A$ topologically finitely generated pro-$l$ groups;
$\phi:H \migisurj A$ a (continuous) surjective homomorphism.
Suppose further that $A$ is abelian, and that $G$ is an $l$-adic Lie group.
%---------------------------------------------------------------------------------------------------------- 
\be
\arrangeVa\item[(i)]
We shall refer to as the \textit{central filtration $\{H(n)\}_{n\geq1}$ on $H$ with respect to
the homorphism $\phi$} the filtration defined as follows: 

\[ H(1):=H \]
\[ H(2):= \mathrm{Ker}(\phi) \] 
\[ H(m):=\bigl\langle [ H(m_1),H(m_2) ] \ \big| \ m_1+m_2=m \bigr\rangle \ \mathrm{for} \ m \geq 3 \]
 --- where $\langle N_i \ | \ i \in I \rangle$ is the group topologically generated by the
$N_i$'s ---.
\\In the following, 
for $a,b,n \in \mbZ$ such that $1 \leq a \leq b, n \geq 1$, we shall write
\[ H(a/b):=H(a)/H(b) \]
\[ \mr{Gr}(H) := \bigoplus_{m \geq 1}H(m/m+1) \]
\[ \mr{Gr}(H)(a/b) := \bigoplus_{b > m\geq a}H(m/m+1)  \]
\[ \mr{Gr}_{\mbQ_l}(H) := \mr{Gr}(H) \otimes_{\mbZ_l} \mbQ_l \]
\[ \mr{Gr}_{\mbQ_l}(a/b) := \mr{Gr}(H)(a/b) \otimes_{\mbZ_l} \mbQ_l \]
\[ H(a/\infty):= \varprojlim_{b > a} H(a/b) \ . \]
%-----------------------------------------------------------------------------------------------
\arrangeVa\item[(ii)]
We shall denote by Lie($G$) the Lie algebra over $\mbQ_l$ determined by the $l$-adic Lie group $G$.
We shall say that $G$ is \textit{nilpotent} if there exists a positive integer $m$ such that 
if we denote by $\{ G(n) \}$ the central filtration with respect to the natural surjection $G \migisurj G^{\mr{ab}}$ (cf.\,(i)), then $G(m)$=\{1\}.
If $G$ is nilpotent, then Lie($G$) is a nilpotent Lie algebra over $\mbQ_l$, hence
determines a connected, unipotent linear algebraic group Lin($G$), which we shall refer to as
the \textit{linear algebraic group associated to} $G$.
In this situation, there exists a natural (continuous) homomorphism (with open image)
\[ G \longmigi \mr{Lin}(G) (\mbQ_l) \]
(from $G$ to the $l$-adic Lie group determined by the $\mbQ_l$-valued points of Lin($G$))
which is uniquely determined (since Lin($G$) is connected and unipotent) by the condition that 
it induce the identity morphism on the associated Lie algebras.
\\ \ \ In the situation of (i), if $1 \leq a \in \mbZ$, then we shall write
\[ \mr{Lie}(H(a/\infty)):= \varprojlim_{b > a} \mr{Lie}(H(a/b)) \]
\[ \mr{Lin}(H(a/\infty)):= \varprojlim_{b > a} \mr{Lin}(H(a/b)) \]
 --- where we note that each $H(a/b)$ is a nilpotent $l$-adic Lie group ---.
\ee
\ede
%\vspace{2mm}
%-----------------------------------------------------------------[end definition]------------------------
%--------------------------------------------------------------[begin definition]------------------------
\bde
\leavevmode\\
 \ \ \ For $n \in \mbZ_{\geq 1}$, we shall denote by 
\[ \{\Delta_{X_n}(m) \} \]
the central filtration of $\Delta_{X_n}$ with respect to the natural surjection
$ \Delta_{X_n} \migisurj \Delta^{\mr{ab}}_{\overline{X}^{\times n}} $ (where $\overline{X}$ denotes the smooth compactification of $X$\,(cf. \S\,0)), and refer to it as the \textit{weight filtration} on $\Delta_{X_n}$.
\ede
%--------------------------------------------------------------[end definition]------------------------- 
%\vspace{1mm}
%--------------------------------------------------------------[begin proposition]------------------
\bpr\leavevmode\\
%(i)
\ \ \ If we equip $\Delta^k_{X_{n/n-1}}$ with the central filtration induced from the identification given by Proposition 1.4 (i) and its weight filtration, then the sequence of morphisms of graded Lie algebras
\[ 1 \longmigi \mr{Gr}(\Delta^k_{X_{n/n-1}}) 
\stackrel{\mr{Gr}(i_k^\Delta)}{\longmigi} \mr{Gr}(\Delta_{X_n}) \stackrel{\mr{Gr}(p_k^\Delta)}{\longmigi} \mr{Gr}(\Delta_{X_{n-1}}) \longmigi 1 \]
induced by the second displayed exact sequence of Definition 1.3 is exact. 
\epr
%---------------------------------------------------------------------[begin proof]-------------------
%\vspace{2mm}
\bpf
See ~\cite{HSHcusp}, Proposition 4.1.
\epf
%--------------------------------------------------------------------[end proposition]------------------
%\vspace{4mm}
Next, let us fix a section $\sigma : G_K \migi \Pi_{X_n} $ of the surjection $ \Pi_{X_n} \migisurj G_K $ arising from the structure morphism of $X_n$.  This section $\sigma$
determines the action of $G_K$ on $\Delta_{X_n}$ by conjugation, 
hence also on
\[ \mr{Gr}_{\mbQ_l}(\Delta_{X_n} )(a/b) , \ \ \mr{Lie}(\Delta_{X_n}(a/b)), \ \ \mr{Lin}(\Delta_{X_n}(a/b))(\mbQ_l ) \]
where $a, b \in \mbZ $ such that $1 \leq a \leq b$.
%\vspace{4mm}
%---------------------------------------------------------------[begin proposition]-------------------
%\vspace{2mm}
\bpr
\leavevmode\\
 \ \ \ Let us assume that $K$ is a finite field whose cardinality we denote by $q_K$, and write $Fr \in G_K$ for the Frobenius element of $G_K$.  Then, relative to the natural conjugate actions determined by $\sigma$:
\vspace{0mm}
\be
\arrangeVa\item[(i)]
The eigenvalues of the action of $Fr$ on $\mr{Lie}_{X_n}(a/a+1)$ are algebraic numbers all of whose complex absolute values are equal to $q^{a/2}_K$ (i.e., weight a).
\arrangeVa\item[(ii)]
There is a unique $G_K$-equivariant isomorphism of Lie algebras
\[ \mr{Lie}(\Delta_{X_n}(a/b)) \isom \mr{Gr}_{\mbQ_l}(\Delta_{X_n})(a/b) \]
which induces the identity isomorphism 
\[ \mr{Lie}(\Delta_{X_n}(c/c+1)) \isom \mr{Gr}_{\mbQ_l}(\Delta_{X_n})(c/c+1) \]
for all $c \in \mbZ_{\geq 1}$
such that $a \leq c < b$.
\ee
\epr
%-------------------------------------------------------------------------[begin proof]------------------
%\vspace{2mm}
\bpf Assertion (i) follows from the ``Riemann hypothesis for abelian varieties over finite fields" (cf., e.g., ~\cite{MUM}, p. 206). Assertion (ii) follows formally from assertion (i) by considering the eigenspaces with respect to the action of $Fr$.
\epf
%--------------------------------------------------------------------[end proposition]-------------------
%\vspace{4mm}
The following proposition is a special case of a result proven previously (cf. ~\cite{STAB}). 
For simplicity, we discuss only the case used in the proofs of the present paper.
%------------------------------------------------------------------------[begin proposition]------------
%\vspace{3mm}
\bpr\leavevmode\\
 \ \ \ For $n=1,2$, the graded Lie algebra $\mr{Gr}( \Delta_{X_n})$  has the following presentation.
\be
\vspace{0mm}
\arrangeVa\item[(i)]
The case $n=1$ (i.e., \,$X_n = X$): \vspace{3mm} 
%--------------------------------------------------------------------------------------------------------
\begin{itemize} \item[] generators ($1 \leq j \leq r, \ 1 \leq i \leq g$) \vspace{3mm}
\begin{itemize}
% generators :  \hspace{5mm}
 \item[$\bullet_1$]  $\zeta_j \ \in \Delta_X(2/3)$ \vspace{3mm}%\notag  \\[9pt]
 \item[$\bullet_2$]  $\alpha_i, \ \beta_i \in \Delta_X(1/2)$ \vspace{6mm}%\notag 
 \end{itemize}
%---------------------------------------------------------------------------------------------------------
\item[] relation \vspace{3mm}
\begin{itemize}
\item[$\bullet_1$] $\sum_{j = 1}^r \zeta_j + \sum_{i = 1}^g [\alpha_i, \beta_i] = 0$ \vspace{6mm} %\notag 
\end{itemize}
\end{itemize}
%-----------------------------------------------------------------------------------------------------------
 --- where $\zeta_j \,(j=1,2,\cdots ,r)$ topologically generates the inertia subgroup in $\Delta_{X}$
(well-defined up to conjugacy) associated to the $j$-th cusp [relative to some ordering of the cusps of $X \times_K \overline{K}$].
\arrangeVa\item[(ii)]
The case $n=2$: \vspace{3mm}
%-------------------------------------------------------------------------------------------------------------
\begin{itemize}
\item[] generators ($1 \leq j \leq r, \ 1 \leq i \leq g, \ k = 1, 2$) \vspace{3mm}
\begin{itemize}
\item[$\bullet_1$] $\zeta \in \Delta_{X_2}(2/3)$ \vspace{3mm}%\notag \\[9pt] 
\item[$\bullet_2$] $  \zeta ^k_j \in \Delta^k_{X_{2/1}}(2/3)$ \vspace{3mm} %\notag \\[9pt]
\item[$\bullet_3$] $\alpha^k_i, \,\beta^k_i \in \Delta^k_{X_{2/1}}(1/2)$ \vspace{6mm} %\notag 
\end{itemize}
%------------------------------------------------------------------------------------------------------------
\item[] relations ($1 \leq j , j'\leq r, \ j \neq j', \ 1 \leq i , i'\leq g, \  \{k, k' \} = \{ 1, 2 \}$)  \vspace{3mm}
\begin{itemize}
\item[$\bullet_1$] $\zeta + \sum_{j =1}^r \zeta^k_j + \sum_{i = 1}^g [\alpha^k_i, \beta^k_i] =0$ \vspace{3mm} %\notag \\[5pt]
\item[$\bullet_2$] $ [\alpha_i^k, \zeta_j^{k'}] = [\beta_i^k, \zeta_j^{k'}]= 0$ \vspace{3mm} %\notag \\[9pt]
\item[$\bullet_3$] $ [ \zeta_j^k, \zeta_{j'}^{k'}] =0$
 \vspace{3mm} %\notag \\[9pt] 
\item[$\bullet_4$] $ [\alpha_i^k, \alpha_{i'}^{k'}]= [\beta_i^k, \beta_{i'}^{k'}] = 0$ \vspace{3mm} 
\item[$\bullet_5$] $ [\alpha_i^k, \beta_{i'}^{k'} ] =
 \begin{cases}
 \zeta & \text{( if $i=i'$)} \\
  0 & \text{( if $i \neq i'$)} 
 \end{cases}$
   \vspace{6mm} 
   \end{itemize} \end{itemize}
%------------------------------------------------------------------------------------------------------------
\vspace{-1mm}
 --- where $\zeta$ topologically generates the image in $\Delta_{X_2}(2/3)$ of the inertia subgroup in $\Delta_{X_2}$
(well-defined up to conjugacy) associated to the diagonal divisor of $X \times_K X$, 
and $\zeta^k_j$ generates the image in $\Delta_{X_{2/1}}^k(2/3)$ of the inertia subgroup in $\Delta^k_{X_{2/1}}$ asssociated to the $j$-th cusp [relative to some ordering of the cusps of $X \times_K \overline{K}$] of the $k$-th factor of $X_2$.
%a different cusp from each other.
\ee
\epr
%--------------------------------------------------------------------[end propositon]-----------------
%\bpf This follows from [6], Proposition ?
%\epf
\leavevmode\\
%%%%%%%%%%%%%%%%%%%%%%%%%%%%%%%%%%%%%%%%%%%%%%%%%
%%%%%%%%%%%%%%%%%%%%%%%%%%%%----[begin section2]----%%%%%%%%%
\section{Switching morphism on configuration spaces}\leavevmode\\
 \ \ \ We continue to use the notations of Section 1.  In this section, we shall introduce certain closed subschemes of $\overline{X}^{\mr{log}}_2 $ equipped with induced log structures --- denoted by $\mbD^{\mr{log}}$ and $\overline{X}_x^{\mr{log}}$ --- and consider various automorphisms induced by the automorphism of $\overline{X}^{\mr{log}}_2 $ determined by  switching the two factors of $X$. The geometry of such log schemes allows us to prove the uniqueness of  certain specific conjugates of induced switching morphisms between fundamental groups that satisfy certain conditions.
 This uniqueness (Proposition 2.5) plays a key role in the proof of Theorem A. 
\\[6mm]
%-----------------------------------------------------------------------------------------------------------
 \ \ \ First, we define a log scheme 
 \[ \mbD^{\mr{log}} \]
to be the log scheme obtained by equipping the diagonal divisor $\overline{X} \subseteq \overline{X}_2$ (which is the restriction of the (1-)morphism $\overline{\mcM}_{g,[r]+1} \migi \overline{\mcM}_{g,[r]+2}$ obtained by gluing the tautological family of curves over $\overline{\mcM}^{\mr{log}}_{g,[r]+1}$ to a trivial family of tripods along the final ordered marked section) with the log structure pulled back from $\overline{X}^{\mr{log}}_2 $. 
Thus, if we write $d : \mbD^{\mr{log}} \migi \overline{X}^{\mr{log}}_2$ for the natural diagonal embedding, then it follows immediately from the definitions that $p_1 \circ d = p_2 \circ d : \mbD^{\mr{log}} \migi \overline{X}^{\mr{log}}$ is a morphism of type $\mbN$ (cf. ~\cite{HSHlog}), i.e., the underlying morphism of schemes is an isomorphism, and the relative characteristic sheaf (cf. \S\,0\,) is locally constant with stalk isomorphic to $\mbN$.

 \ Observe that the (1-)automorphism on $\overline{\mcM}^{\mr{log}}_{g,[r]+2}$ over $\overline{\mcM}^{\mr{log}}_{g,[r]}$ given by switching the two ordered marked points of the tautological family of curves over $\overline{\mcM}^{\mr{log}}_{g,[r]+2}$
induces automorphisms $s$, $\overline{s}$, and $s_\mbD$, which fit into a commutative diagram as follows:

\[ \hspace{35mm} \begin{CD}
\mbD^{\mr{log}}
@> d >>
\overline{X}^{\mr{log}}_2
@> p := (p_1,p_2) >> \overline{X}^{\mr{log}} \times_K \overline{X}^{\mr{log}}
\\
@V s_\mbD VV @V s VV @V \overline{s} VV 
\\
\mbD^{\mr{log}}
@> d >>
\overline{X}^{\mr{log}}_2
@> p:=(p_1,p_2)  >>
\overline{X}^{\mr{log}} \times_K \overline{X}^{\mr{log}}  .
\end{CD} \hspace{25mm} (*)^X \]

%----------------------------------------------------------------------------[begin lemma]---------------
\ble \leavevmode\\
 \ \ \ In the notation of the above situation,
\vspace{-1mm}
\be
\arrangeVa\item[(i)] 
$\overline{s}$ is the morphism determined by switching the two factors.
\arrangeVa\item[(ii)]
 $s_\mbD$ is the identity morphism on the underlying scheme; on the sheaf of monoids defining the log structure of $\mbD^{\mr{log}}$, for any \'{e}tale local section $s$ of $\mcM_{\mbD}$ such that $``s =0"$ defines the diagonal divisor $\overline{X} \subseteq \overline{X}_2 $, 
\[ s_\mbD (s) = -s \ . \]
\ee
\ele
%--------------------------------------------------------------------------------[begin proof]--------------
\bpf
Recall (cf. ~\cite{KNUD}, Corollary 2.6) that $\overline{X}_2$ is obtained by blowing-up $\overline{X} \times_K \overline{X}$ along the intersection of the diagonal divisor and the pull-backs of  the cusps via $p_1, p_2 : \overline{X}_2 \migi \overline{X}$.  Thus, one verifies easily that assertions (i) and (ii) follow immediately from the fact  that the ring homomorphism corresponding to $\overline{s}$ in an affine neighborhood of any diagonal point may be expressed as 
\[  A \otimes_K A \longmigi A \otimes_K A  \]
\[ \sum_j a_j \otimes a'_j \mapsto \sum_j a'_j \otimes a_j \ ,\]
hence  maps $s$ to $- s$ for any local section $s$ such that $``s =0"$ defines the diagonal divisor $\overline{X} \subseteq \overline{X} \times_K \overline{X}$. 
\epf
%--------------------------------------------------------------------------------[end lemma]-------------
%---------------------------------------------------------------------------------[begin remark]----------
\begin{rema} \leavevmode\\
 \ \ \ Lemma 2.1 (ii) can be interpreted as the assertion that the automorphism induced by $s_\mbD$ 
on the sheaf of monoids $\mcM_{\mbD}$ defining the log structure of $\mbD^{\mr{log}}$
 may be expressed, relative to  the \'{e}tale local splitting of $\mcM_{\mbD} \migisurj \mcM_{\mbD} / \mcO^*_X \cong \mbN$ corresponding to $s$, as
\[ \mbN \ \oplus \ \mcO_X^* \longisom \mbN \ \oplus \ \mcO_X^* \]
\[ \ \ \ \ \ \ (m, v) \longmapsto (m, (-1)^m v) \ . \]
\end{rema}
%-------------------------------------------------------------------------------[end remark]------------
\vspace{5mm}
 Next, we introduce the log scheme $\overline{X}^{\mr{log}}_x$ that appears in the discussion at the beginning of this section.
 Let $x^{\mr{log}} \migi \overline{X}^{\mr{log}}$ be a strict  morphism  (cf. ~\cite{ILL}, 1.2) such that  the underlying scheme of $x^{\mr{log}}$ is $K$-isomorphic to $\mr{Spec}(K)$.  
We shall write 
\[ \overline{X}_x^{\mr{log}} :=  x^{\mr{log}} \times_{\overline{X}^{\mr{log}}} \overline{X}_2^{\mr{log}}, \]
\[ \tilde{x}^{\mr{log}} := x^{\mr{log}} \times_{\overline{X}^{\mr{log}}} \mbD^{\mr{log}}, \]
--- where the morphism $\overline{X}_2^{\mr{log}} \migi \overline{X}^{\mr{log}}$ (resp., $\mbD^{\mr{log}} \migi \overline{X}^{\mr{log}}$) in the fiber product defining $\overline{X}_x^{\mr{log}}$ (resp., $\tilde{x}^{\mr{log}}$) is $p_1$ (resp., $p_1 \circ d = p_2 \circ d$)
 --- and refer to $\overline{X}_x^{\mr{log}}$ (resp., $\tilde{x}^{\mr{log}}$) as the \textit{cuspidalization of $X$ at $x$} (resp., \textit{diagonal cusp of $\overline{X}_x^{\mr{log}}$}).
We note that both the log structure of $x^{\mr{log}}$ and the underlying scheme of $\overline{X}_x^{\mr{log}}$ depend on the choice of $x \in \overline{X}$:
\be
\arrangeVa\item[(1)] 
{\it The Case $x \in X$}: \\
In this case, $x = x^{\mr{log}}$, i.e., the log structure of $x^{\mr{log}}$ is trivial. 
As we discussed in Section 1, the underlying scheme of $\overline{X}_x^{\mr{log}}$ is naturally isomorphic to $\overline{X}$; this isomorphism maps $\tilde{x}$ to $x$ and the interior of $\overline{X}_x^{\mr{log}}$ onto $X \setminus \{ x \}$. 
\arrangeVa\item[(2)]
{\it The Case $x \in \overline{X} \setminus X$}: \\
In this case, the log structure of $x^{\mr{log}}$ has a chart modeled on $\mbN$, which determines a local uniformizer of $X$ at $x$.   The scheme
$\overline{X}_x$ consists of precisely two irreducible components, one of which maps to the point $x \in X$ (resp., maps isomorphically to $\overline{X}$) via $\overline{X}_x^{\mr{log}} \stackrel{p_2 \circ i_1}{\longmigi} \overline{X}^{\mr{log}}$; denote this irreducible component by $\overline{\mbP}_K$ (resp., $\overline{X}$, via a slight abuse of notation).  Thus, $\overline{X}$, $\overline{\mbP}_K$ are joined at a single \textit{node} $\nu_x$. Let us refer to $\overline{X}$ (resp., $\overline{\mbP}_K$, $\nu_x$) as the \textit{major cuspidal component} (resp., the \textit{minor cuspidal component}, the \textit{nexus}) at $x$, and denote by
$\overline{X}^{\mr{log}'}$, $\overline{\mbP}_K^{\mr{log}'}$, $\nu_x^{\mr{log}}$ the log schemes obtained by equipping 
$\overline{X}$, $\overline{\mbP}_K$, $\nu_x$ with the respective log structures pulled back from $\overline{X}_x^{\mr{log}}$ (cf. ~\cite{MZKcomb2}, Definition 1.4).
Note that the 1-interior of $\overline{X}^{\mr{log}'}$ (resp., $\overline{\mbP}_K^{\mr{log}'}$) is naturally isomorphic to $X$ (resp., is a \textit{tripod}).
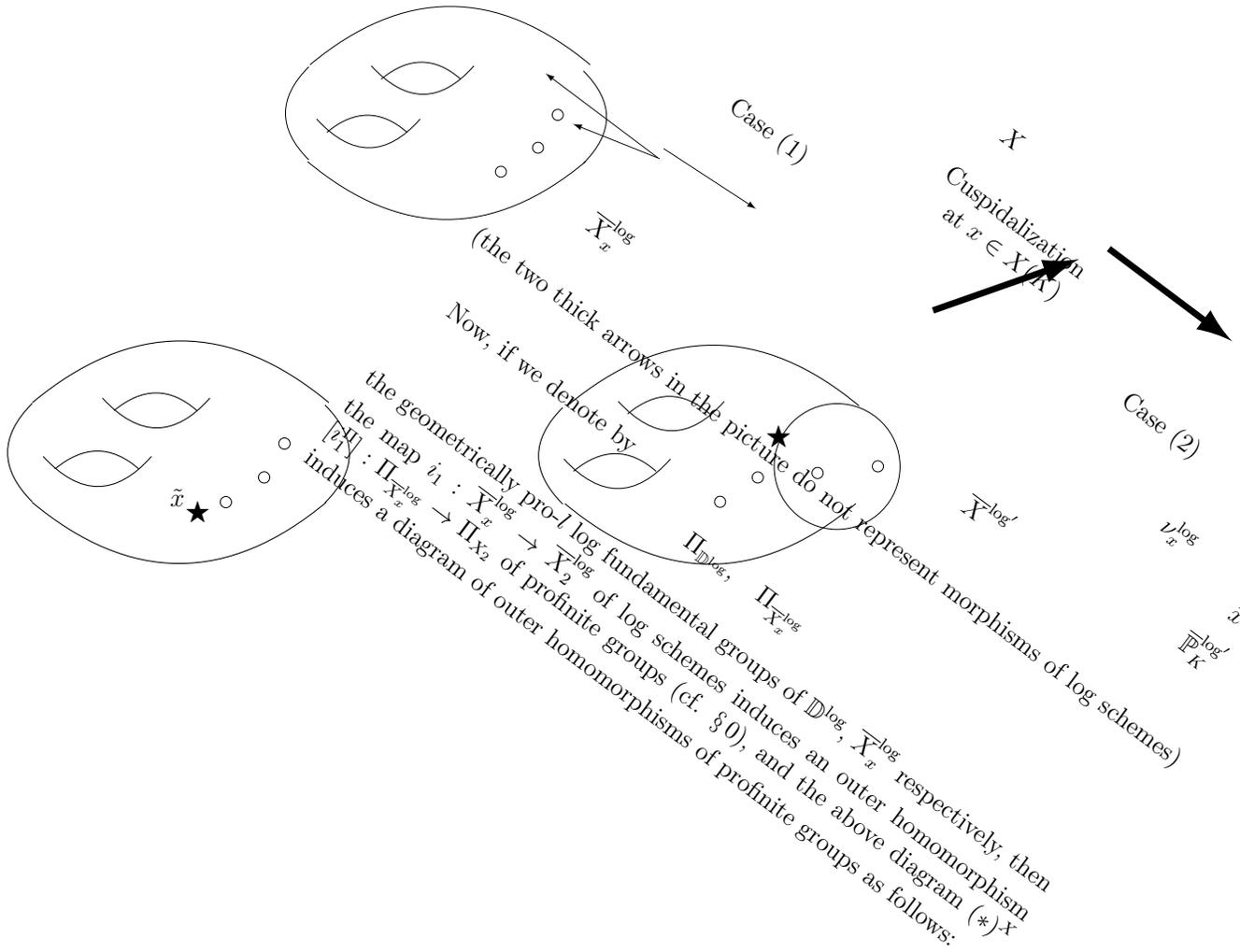
\ee
\begin{center}
\begin{picture}(400,270)
\qbezier(20,90)(80,140)(140,90)
\qbezier(20,90)(0,70)(20,50)
\qbezier(20,50)(80,0)(140,50)
\qbezier(140,50)(160,70)(140,90)

\qbezier(360,90)(300,140)(240,90)
%\qbezier(360,90)(380,70)(360,50)
\qbezier(360,50)(300,0)(240,50)
\qbezier(240,50)(220,70)(240,90)

\qbezier(25,70)(48,45)(70,70)
\qbezier(30,65)(48,80)(65,65)
\qbezier(50,95)(73,70)(95,95)
\qbezier(55,90)(73,105)(90,90)

\qbezier(245,70)(268,45)(290,70)
\qbezier(250,65)(268,80)(285,65)
\qbezier(270,95)(293,70)(315,95)
\qbezier(275,90)(293,105)(310,90)

\qbezier(140,240)(200,290)(260,240)
\qbezier(140,240)(120,220)(140,200)
\qbezier(140,200)(200,150)(260,200)
\qbezier(260,200)(280,220)(260,240)

\qbezier(145,220)(168,195)(190,220)
\qbezier(150,215)(168,230)(185,215)
\qbezier(170,245)(193,220)(215,245)
\qbezier(175,240)(193,255)(210,240)

\put(133,80){\circle{5}}
\put(124,65){\circle{5}}
\put(107,54){\circle{5}}
%\put(125,85){x}

\put(348,80){$\bigstar$}
\put(371,67){\circle{5}}
\put(398,70){\circle{5}}

\put(344,65){\circle{5}}
\put(327,54){\circle{5}}

\put(380,70){\circle{57}}

\put(253,230){\circle{5}}
\put(244,215){\circle{5}}
\put(227,204){\circle{5}}

\put(82,44){$\bigstar$}
\put(75,50){$\tilde{x}$}
\put(300,210){\vector(-5,2){40}}
\put(300,205){\vector(-5,1){50}}
\put(300,200){\vector(-15,1){65}}
\linethickness{3pt}

\put(155,180){\vector(-2,-3){40}}
\put(242,180){\vector(2,-3){40}}
\put(160,160){Cuspidalization}
\put(170,145){at $x \in X(K)$}
\put(340,89){$\nu_x^{\mr{log}}$}
\put(60,125){Case (1)}
\put(290,125){Case (2)}
\put(390,75){$\tilde{x}$}
\put(170,190){$X$}
\put(40,40){$\overline{X}_x^{\mr{log}}$}
\put(260,40){$\overline{X}^{\mr{log}'}$}
\put(380,50){$\overline{\mbP}_K^{\mr{log}'}$}
\put(305,202){cusps}
\end{picture}
\end{center}
 \vspace{-10mm} \[ \text{(the two thick arrows in the picture do not represent morphisms of log
schemes)} \]
\vspace{0mm}
%----------------------------------------------------------------------------------------------------------

Now, if we denote by 
 \[ \Pi_{\mbD^{\mr{log}}}, \ \ \ \Pi_{\overline{X}^{\mr{log}}_x} \]
   the geometrically pro-$l$ log fundamental groups of $\mbD^{\mr{log}}$, $\overline{X}^{\mr{log}}_x$ respectively, then the map $i_1:\overline{X}_x^{\mr{log}} \migi \overline{X}_2^{\mr{log}}$ of log schemes induces an outer homomorphism $[i^\Pi_1 ]: \Pi_{\overline{X}^{\mr{log}}_x} \migi \Pi_{X_2}$ of profinite groups (cf. \S\,0),
%(i.e., an element of $\mr{Hom}(\Pi_{\overline{X}^{\mr{log}}_x}, \Pi_{X_2}) /\sim$, where $\sim$ denotes the relation determined by  $f \sim \mr{Inn}(s) \circ f$ for some $s \in \Pi_{X_2}$),
and the above diagram $(*)^X$ induces a diagram of outer homomorphisms of profinite groups as follows:
\\
\[ \hspace{35mm} \begin{CD}
\Pi_{\mbD^{\mr{log}}} @> [ d^\Pi ] >> \Pi_{X_2} @> [ p^\Pi ] >> \Pi_X \times_{G_K} \Pi_X
\\
@V [ s_\mbD^\Pi ] V \wr V @V [ s^\Pi ] V \wr V @V [ \overline{s}^\Pi ] V \wr V
\\
\Pi_{\mbD^{\mr{log}}} @> [ d^\Pi ] >> \Pi_{X_2} @> [ p^\Pi ] >> \Pi_X \times_{G_K} \Pi_X \ . 
\end{CD} \hspace{25mm} (*)^\Pi \] 
\\
Note that the homomorphisms corresponding to the arrow $[i^\Pi_1]$ and the arrows in the diagram $(*)^\Pi$ are only defined (i.e., in the absence of appropriate choices of basepoints of respective log schemes) up to conjugacy, and that $[\overline{s}^\Pi]$ coincides
with the morphism obtained by switching the two factors.
The main purpose of this section is to give characterizations of certain specific
choices within these conjugacy classes of homomorphisms.
%----------------------------------------------------------------------[begin definition]-----------------
%\vspace{4mm}
\bde \leavevmode\\
\be \vspace{-5mm} \arrangeVa\item[(i)] 
We shall denote by 
 \[ (C_{i_1}^X)  \  \bigl( \text{resp., } (C_{p_1}^X), (C_{p_2}^X), (C_{d}^X)\bigr) \]
 a {\it choice of a specific homomorphism}
[i.e., in the sense that it is {\it not} subject to conjugacy indeterminacy!] 
 \[  i_1^\Pi: \Pi_{\overline{X}_x^{\mr{log}}} \longmigi \Pi_{X_2} \]
\[  \bigl(\text{resp.,}  \  p_1^\Pi: \Pi_{X_2} \longmigi \Pi_X, \
 p_2^\Pi: \Pi_{X_2} \longmigi \Pi_X, \
 d^\Pi: \Pi_{\mbD^{\mr{log}}} \longmigi \Pi_{X_2} \bigr) \]
induced by the morphism of log schemes
$i_1:\overline{X}_x^{\mr{log}} \migi \overline{X}_2^{\mr{log}}$ (resp., $p_1:\overline{X}_2^{\mr{log}} \migi \overline{X}^{\mr{log}}$,  $p_1:\overline{X}_2^{\mr{log}} \migi \overline{X}^{\mr{log}}$, $\mbD^{\mr{log}} \migi \overline{X}_2^{\mr{log}}$). 
\arrangeVa\item[(ii)] 
We shall denote by 
\[ (C_{\tilde{x}}^X) \ \bigl( \text{resp., } (C_D^X) \big)\]
 a {\it choice of a specific subgroup} --- i.e., of a {\it specific decomposition group} ---
\[  D_{\tilde{x}} \subseteq \Pi_{\overline{X}_x^{\mr{log}}} \ \bigl( \text{resp., } D_X \subseteq \Pi_{X_2} \bigr) \]
associated to $\tilde{x}^{\mr{log}}$ of $\overline{X}_x^{\mr{log}}$ 
(resp., the diagonal divisor of $\overline{X}_2$),
 among the various conjugates of this subgroup.  
 Note that such a choice determines  a {\it choice of a specific subgroup}
 --- i.e., of a {\it specific inertia group} ---
\[ I_{\tilde{x}} := D_{\tilde{x}} \cap \Delta_{X_{2/1}^1} \subseteq \Pi_{\overline{X}_x^{\mr{log}}} \ \bigl( \text{resp., } I_X := D_X \cap \Delta_{X_{2/1}^1} \subseteq \Pi_{X_2} \bigr) \]
among the various conjugates of this subgroup.
\arrangeVa\item[(iii)] 
 Assume that we have fixed a choice $(C_{\tilde{x}}^X)$ of a specific decomposition group $D_{\tilde{x}} \subseteq \Pi_{\overline{X}_x^{\mr{log}}}$ (hence also of a specific inertia group $I_{\tilde{x}} \subseteq \Pi_{\overline{X}_x^{\mr{log}}}$).
Then we shall denote by  
\[ (C_\sigma^{(-) X}) \]
a \textit{choice of a specific section}
\[ \hspace{30mm} \sigma^{(-)} : G_K^{(-)} \longmigi D_{\tilde{x}}^{(-)} \hspace{40mm} \]
--- where the symbol $(-)$ denotes either the presence or absence of $``\dagger"$
[thus,  a choice $(C_\sigma^X)$ determines a unique choice $(C_\sigma^{\dagger X})$ by 
restriction]  --- of the natural surjection $D_{\tilde{x}}^{(-)} \migi G_K^{(-)}$ (cf. Remark 2.2.1)
and by 
\[ (C_\delta^X) \]
a \textit{choice of a specific 1-cocycle map}
\[ \hspace{30mm} \delta : G^\dagger_K  \longmigi I_{\tilde{x}} \hspace{40mm} \]
 representing the Kummer class $-1 \in (K^\times )^\wedge$ (cf. Remark 2.2.1).
\ee
\ede
%-------------------------------------------------------------------------[end definition]---------------
%\vspace{3mm}
Before proceeding, we pause to make a remark concerning Definition 2.2.
%-------------------------------------------------------------------------[begin remark]---------------
\begin{rema}
\leavevmode\\
\be \vspace{-4mm}
\arrangeVa\item[(i)] Recall that the natural surjection $D_{\tilde{x}} \migisurj G_K$ (which, 
since $G_K$ is {\it abelian}, is {\it uniquely determined} without any conjugacy indeterminacies) has a section. Indeed, when $x \in X$ (resp., $x \in \overline{X} \setminus X$), fixing a choice of such a  section is equivalent to extracting roots of any local uniformizer(s) of the divisor(s) $\mbD \subseteq \overline{X}_2$ (resp., $\mbD \subseteq \overline{X}_2$ and $\overline{X}_x \subseteq \overline{X}_2$) at $\tilde{x}$. 
\arrangeVa\item[(ii)] 
We shall consider the restriction map
$H^1(G_K, I_{\tilde{x}}) \migi H^1(G_K^\dagger, I_{\tilde{x}})$
of cohomology groups induced by the natural inclusion
$G_K^\dagger \migiincl G_K$.
Since $G_K^\dagger$ is the maximal pro-$l$ subgroup of $G_K$
and $I_{\tilde{x}}$ is isomorphic to $\mbZ_l(1)$ as a $G_K$-module,
this restriction map determines an isomorphism of $H^1(G_K, I_{\tilde{x}})$ 
with $H^1(G_K^\dagger, I_{\tilde{x}})$,
hence also with the maximal pro-$l$ completion $(K^\times )^{\wedge}$ of the multiplicative group $K^\times$ of $K$. 
Therefore, if we denote by $Z^1(G_K^\dagger, I_{\tilde{x}})$ (resp., $Z^1(G_K, I_{\tilde{x}})$) the set of (continuous) 1-cocycles  of $G_K^\dagger$ (resp., $G_K$) with coefficients in $I_{\tilde{x}}$, 
then it makes sense to refer to any element of $Z^1(G_K^\dagger, I_{\tilde{x}})$ (resp., $Z^1(G_K, I_{\tilde{x}})$) belonging to the inverse image of $a \in (K^\times )^{\wedge} \cong H^1(G_K^\dagger, I_{\tilde{x}})$ (resp., $\cong H^1(G_K^\dagger, I_{\tilde{x}}))$, via the natural surjection, as a \textit{(continuous) 1-cocycle representing the Kummer class $a$}.
\ee
 \end{rema}
%-------------------------------------------------------------------------------[end remark]-------------

%----------------------------------------------------------------------------[begin lemma]--------------
%\vspace{3mm}
\ble \leavevmode\\
%\vspace{-5mm}
 \ \ \ For any choice $(C_{\tilde{x}}^X)$ (resp., $(C_D^X)$) of a specific decomposition group $D_{\tilde{x}} \subseteq \Pi_{\overline{X}_x^{\mr{log}}}$ (resp., $D_X \subseteq \Pi_{X_2}$),
$I_{\tilde{x}}$ (resp., $I_X$) is normally terminal in $\Delta_{X_{2/1}}^1$ (cf. \S0),
and $D_{\tilde{x}}^{(-)}$ (resp., $D_X^{(-)}$) --- where the symbol $(-)$ denotes either the presence or absence of $``\dagger"$ --- coincides with $N_{\Pi_{\overline{X}_x^{\mr{log}}}}(I_{\tilde{x}})^{(-)}$ (resp., $N_{\Pi_{X_2}}(I_X)^{(-)}$)(cf. \S0).
 \ele
 %-------------------------------------------------------------------------[begin proof]------------- 
\bpf
Recall that, by definition, we have $I_{\tilde{x}} = D_{\tilde{x}} \cap \Delta_{X_{2/1}}^1 \subseteq \Pi_{X_2}$ and $I_X = D_X \cap \Delta_{X_{2/1}}^1 \subseteq \Pi_{X_2}$.
Next, let us recall the well-known fact (cf., e.g., ~\cite{NAK}, (2.3.1)) that $I_{\tilde{x}}$ and $I_X$ are normally terminal (cf. \S0) in $\Delta_{X_{2/1}}^1$. 
Thus, the \textit{resp'd} assertion follows immediately from the fact that $p^\Pi_1$ maps $D_X$ onto $\Pi_X$.
On the other hand, the \textit{non-resp'd} assertion follows immediately from the observation that the images of $D_{\tilde{x}}$ and $\Pi_{\overline{X}_x^{\mr{log}}}$ via $p_1^\Pi \circ i_1^\Pi$ {\it coincide}.  This observation is a consequence of the geometry
of the corresponding morphisms of log schemes, which implies that both of these images coincide with a decomposition group $\subseteq \Pi_X$ associated to the point $x$.
\epf
%-----------------------------------------------------------------------[end lemma]----------------------
%-------------------------------------------------------------------------[begin lemma]-------------------
\ble \leavevmode\\
\be \vspace{-5mm}
\arrangeVa\item[(i)]
If we fix a choice $(C_d^X)$ of $d^\Pi : \Pi_{\mbD^{\mr{log}}} \migi \Pi_{X_2}$, then there exists a unique choice $(C_D^X)$ of $D_X \subseteq \Pi_{X_2}$ such that the image of $d^\Pi$ coincides with $D_X$.
By contrast, if we fix a choice $(C_D^X)$ of $D_X \subseteq \Pi_{X_2}$, then there exists a (not necessarily unique!) choice $(C_d^X)$ of $d^\Pi : \Pi_{\mbD^{\mr{log}}} \migi \Pi_{X_2}$ such that the image of $d^\Pi$ coincides with $D_X$.
\arrangeVa\item[(ii)]
If we fix  a triple of choices $(C_{i_1}^X)$,  $(C_{p_1}^X)$ and $(C_{\tilde{x}}^X)$, then there exists a unique pair consisting of 
a choice $(C_{p_2}^X)$ of  $p^\Pi_2 : \Pi_{X_2} \migi \Pi_X$ and a choice $(C_D^X)$ of $D_X \subseteq \Pi_{X_2}$
that satisfy the following conditions:  
\be
\arrangeVa\item[(1)] The image of the inertia group $I_X \subseteq D_X$ in $\Pi_{X_2}$ coincides with the image of $I_{\tilde{x}}$ via $i_1^\Pi$.
\arrangeVa\item[(2)]
The homomorphism $(p^\Pi_1, p^\Pi_2) : \Pi_{X_2} \migi \Pi_X \times_{G_K} \Pi_X$ maps $D_X$ onto the image of the diagonal embedding $\Pi_X\migiincl \Pi_X\times_{G_K} \Pi_X$.
\ee
\ee
\ele
%-------------------------------------------------------------------------[begin proof]--------------
\bpf
Assertion (i) follows immediately from the definitions of $\Pi_{\mbD^{\mr{log}}}$ and $D_X$.
Next, we consider assertion (ii). 
First, let us observe that it follows immediately from the various definitions involved
that $I_X$ and $I_{\tilde{x}}$ are $\Pi_{X_2}$-conjugate. 
Since, by Lemma 2.3, $D_X$ coincides with the normalizer of $I_X$ in $\Pi_{X_2}$,
it suffices to take $D_X$ to be the normalizer of $I_{\tilde{x}}$ in $\Pi_{X_2}$ and $p_2^\Pi$ 
to be such that the condition  $p_1^\Pi = p_2^\Pi $ is satisfied on $D_X$. Uniqueness follows immediately from the conditions (1), (2) and the surjectivity of the restriction of $p^\Pi_2$ to $D_X$.
\epf
%--------------------------------------------------------------------------[end lemma]---------------

%-------------------------------------------------------------------[begin proposition]-----------
%\vspace{4mm}
\bpr\leavevmode\\
 \ \ \ If we fix arbitrary choices $(C_{i_1}^X)$, $(C_{p_1}^X)$, $(C_{\tilde{x}}^X)$, $(C_{\sigma}^{\dagger X})$, and $(C_{\delta}^X)$,
then
there exists a unique triple of choices  
consisting of $(C_{p_2}^X)$, $(C_D^X)$, and a 
{\bf choice of a specific automorphism}  induced by  $s : \overline{X}^{\mr{log}}_2 \isom \overline{X}^{\mr{log}}_2$ 
\[ s^\dagger : \Pi_{X_2}^\dagger \stackrel{\sim}{\longmigi} \Pi_{X_2}^\dagger \] 
 --- which we shall denote by $(C_s^{\dagger X})$ --- satisfying the two conditions (1), (2) stated in Lemma 2.4, (ii), as well as the following conditions:

\begin{itemize}
\item[(1)] 
the morphism 
$\overline{s}^\dagger :\Pi_X^\dagger \times_{G_K^\dagger} \Pi_X^\dagger \longisom \Pi_X^\dagger \times_{G_K^\dagger} \Pi_X^\dagger$ induced by passing to the quotient
$\Pi_{X_2}^\dagger \stackrel{p^\dagger}{\migisurj} \Pi_X^\dagger \times_{G_K^\dagger} \Pi_X^\dagger$ determined by $p_1^\Pi, p_2^\Pi$ coincides with   the morphism obtained by switching the two factors.
\\
\item[(2)]
$s^\dagger$ preserves $D_X^\dagger \subseteq \Pi_{X_2}^\dagger$, and the restriction $s^\dagger |_{D_X^\dagger} : D_X^\dagger \isom D_X^\dagger$ corresponds to an automorphism induced by $s_\mbD : \mbD^{\mr{log}} \isom \mbD^{\mr{log}}$ via the identification $\Pi_{\mbD^{\mr{log}}}^\dagger \isom D_X^\dagger$ determined by some choice of a specific homomorphism $d^\Pi : \Pi_{\mbD^{\mr{log}}} \migi \Pi_{X_2}$ whose image coincides with $D_X$ (cf. Lemma 2.4, (i)).
\\
\item[(3)] 
The continuous function
$G_K^\dagger \migi \Pi_{X_2 }^\dagger$ defined by 
\[ g \mapsto (s^\dagger \circ \sigma^\dagger )(g) \cdot \sigma^\dagger (g)^{-1} \]
is valued in $I_{\tilde{x}} \subseteq \Pi_{X_2 }^\dagger$ and coincides with the ``$\delta$'' determined by $(C_\delta^X)$.
\end{itemize}
In particular, $s^\dagger$ induces the identity morphism on $I_X \subseteq \Pi_{X_2}^\dagger$.
\epr
%-------------------------------------------------------------------[begin proof]----------------
\bpf
We begin by proving the existence portion.
Let us consider the following (not necessarily commutative) diagram
\[ \hspace{35mm} \begin{CD}
\Pi_{\mbD^{\mr{log}}}^\dagger @> d^\dagger >> \Pi_{X_2}^\dagger @> p^\dagger >> \Pi_X^\dagger \times_{G_K^\dagger} \Pi_X^\dagger
\\
@V \underline{s}^\dagger VV @V s^\dagger VV @V \overline{s}^\dagger VV
\\
\Pi_{\mbD^{\mr{log}}}^\dagger @> d^\dagger >> \Pi_{X_2}^\dagger @> p^\dagger >> \Pi_X^\dagger \times_{G_K^\dagger} \Pi_X^\dagger
\end{CD} \hspace{25mm} (*)^\dagger \]
induced by $(*)^\Pi$ consisting of the horizontal arrows arising from the choice $(C_{p_1}^X)$ fixed in advance and the pair of choices $(C_{p_2}^X)$, $(C_d^X)$ obtained by applying Lemma 2.4 (i), (ii), and arbitrary choices of the vertical arrows. 
By the surjectivity of $p^\dagger$, we can take $s^\dagger, \overline{s}^\dagger$ such that the right-hand square of the diagram $(*)^\dagger$
commutes, and condition (1) is satisfied.
The commutativity of the rectangle in $(*)^\dagger$ up to conjugacy
implies that 
there exists $\lambda \in \Pi_X^\dagger \times_{G_K^\dagger} \Pi_X^\dagger$ 
such that $\overline{s}^\dagger \circ (p^\dagger \circ d^\dagger) = \mr{Inn} ( \lambda ) \circ (p^\dagger \circ d^\dagger) \circ \underline{s}^\dagger$
(where $\mr{Inn} ( \lambda )$ denotes the inner automorphism obtained by conjugating by $\lambda$). 
By the construction of the choice $(C_D^X)$ (cf. condition (2) of Lemma 2.4 (ii)), 
$p^\dagger \circ d^\dagger$ maps $\Pi_{\mbD^{\mr{log}}}^\dagger$
onto the  subgroup of diagonal elements of $\Pi_X^\dagger \times_{G_K^\dagger} \Pi_X^\dagger$; thus, $\mr{Inn} ( \lambda )$ preserves this diagonal subgroup. 
Since $\Pi_X^\dagger$ is center-free (by Proposition 1.4 (iii)),
we thus conclude that $\lambda$ is a diagonal element.
Thus, by taking a lifting $\tilde{\lambda} \in \Pi_{\mbD^{\mr{log}}}^\dagger$ of $\lambda$
and replacing $\underline{s}^\dagger$ by $\mr{Inn} ( \tilde{\lambda}^{-1} ) \circ \underline{s}^\dagger$,
we can make the rectangle in  $(*)^\dagger$ commute in the strict sense.
Next, we observe (by applying again the commutativity of the rectangle in $(*)^\dagger$ up to conjugacy) that
$s^\dagger \circ d^\dagger = \mr{Inn} ( \mu ) \circ d^\dagger \circ \underline{s}^\dagger$ for some $\mu \in \Pi_{X_2}^\dagger$.
By the commutativity of the rectangle in $(*)^\dagger$,
$\mu$ projects via $p^\dagger$ into the center of  $\Pi_X^\dagger \times_{G_K^\dagger} \Pi_X^\dagger$,
hence (by Proposition 1.4, (iii)), to the unit element. 
Therefore, by replacing $s^\dagger$ by $\mr{Inn}( \mu^{-1} )  \circ s^\dagger$,
we conclude that we may choose $\underline{s}^\dagger$, $s^\dagger$, and
$\overline{s}^\dagger$ so that the diagram $(*)^\dagger$ {\it commutes}, and, moreover,
conditions (1) and (2) are satisfied. 

Next, observe that by restricting $s^\dagger$ to $D_X^\dagger$, we obtain
a commutative diagram
\bcd
 1 @>>> I_X @>>> D_X^\dagger @> p^\dagger|_{D^\dagger_X} >> \Pi_X^\dagger @>>> 1 
\\
@. @V s^\dagger |_{I_X}  V \wr V @V s^\dagger |_{D_X^\dagger} V \wr V @V\text{id} V \wr V @.
\\
 1 @>>> I_X @>>> D_X^\dagger @> p^\dagger|_{D^\dagger_X} >> \Pi_X^\dagger @>>> 1 
\ecd  
in which the right-hand vertical arrow is the identity automorphism of $\Pi_X^\dagger$.
Write $\mbM \subseteq \mbQ$ for the monoid of positive rational numbers with 
$l$-power denominators,
and $\mcN$ for the monoid of global sections  of the sheaf of monoids defining the log structure on a universal geometrically pro-$l$ k\'{e}t covering (cf. ~\cite{ILL}, Definition 3.1) of $x^{\mr{log}} \times_{\overline{X}^{\mr{log}}} \mbD^{\mr{log}}$.
When $x \in X$ (resp., $x \in \overline{X} \setminus X$), $\mcN$ admits a direct sum decomposition $\mcN \cong \mbM \oplus \overline{K}^\times$ (resp., $\mcN \cong \mbM \oplus \mbM \oplus \overline{K}^\times$), where (cf. Remark 2.2.1 (i)) the first  factor (resp., first two factors) of the direct sum arise(s) from extracting roots of a local uniformizer of the divisor $\mbD \subseteq \overline{X}_2$ at $\tilde{x}$ (resp., of local uniformizers of the two irreducible divisors defining the log structure of $\overline{X}^{\mr{log}}_2$ at $\tilde{x}$) 
in a fashion compatible with the choice  $(C^{\dagger X}_\sigma)$ of $\sigma$.  
Here, in the resp'd case, we assume that the first factor ``$\mbM$'' corresponds to the divisor
$\mbD\subseteq\overline{X}_2$.
Next, observe that it follows from Lemma 2.1 (ii), together with the well-known local
structure of $\overline{X}_2$ in a neighborhood of $\tilde{x}$, that 
 the automorphism of $\mcN \cong \mbM \oplus \overline{K}^\times$ (resp., $\mcN \cong \mbM \oplus \mbM \oplus \overline{K}^\times $) induced by
the automorphism $\underline{s}^\dagger$  of $\Pi_{\mbD^{\mr{log}}}^\dagger$ may be  expressed  in the form 
\[ (\frac{a}{l^{m}}, k) \mapsto (\frac{a}{l^m}, (-1)^{\frac{a}{l^m}} \cdot k) \]
\[ \bigl( \text{resp.,} (\frac{a_1}{l^{m_1}}, \frac{a_2}{l^{m_2}}, k) \mapsto (\frac{a_1}{l^{m_1}}, \frac{a_2}{l^{m_2}}, (-1)^{\frac{a_1}{l^{m_1}}} \cdot k) \bigr) \]
for a suitable choice of a projective system $\{ (-1)^{\frac{1}{l^m}} \}_{m \in {\mbZ_{\geq 0}}}$ of $l$-power roots of $-1$.
In particular, we conclude that the restriction $s^\dagger |_{I_X}$
is the identity morphism, and that
the $1$-cocycle $G_K^\dagger\ni g\mapsto ( s^\dagger \circ \sigma^\dagger )(g) \cdot \sigma^\dagger (g)^{-1}$
 is valued in $I_X = I_{\tilde{x}}$ (cf. condition (1) of Lemma 2.4 (ii)).
 Therefore, by replacing $\underline{s}^\dagger$, $s^\dagger$ by their composites with 
a suitable $I_X$-inner automorphism, we may assume that condition (3) is satisfied. This completes the proof of the existence assertion.

 Next we prove the uniqueness portion.
If $s^\dagger_1$, $s^\dagger_2$ are two maps
 that satisfy conditions (1), (2) and (3), then $s_1^\dagger \circ (s_2^{\dagger})^{ -1} = \mr{Inn}( \eta ) \in \mr{Aut}(\Pi_{X_2}^\dagger)$
for some $\eta \in \Pi_{X_2}^\dagger$,
and it follows from condition (2)
that $\mr{Inn}( \eta )$ preserves the subgroup $D_X^\dagger \subseteq \Pi_{X_2}^\dagger$.
Since $D_X^\dagger$ is normally terminal in $\Pi_{X_2}^\dagger$ (cf. Lemma 2.3),  
we thus conclude that $\eta$ is in $D_X^\dagger$.
Moreover, it follows from condition (1) and the fact that $\Pi_X^\dagger$ is center-free
(cf. Proposition 1.4, (iii)), that $\eta$
lies in $\mr{Ker} \bigl( D_X^\dagger \stackrel{p^\dagger |_{D_X^\dagger}}{\migi} \Pi_{X}^\dagger \times_{G_K^\dagger} \Pi_X^\dagger \bigr)$, i.e., $\eta \in I_X$.
On the other hand, since the section $\sigma^\dagger$ acts {\it faithfully} on $I_X $ via the cyclotomic character, condition (3) implies that
$\eta $ is the unit element, i.e., that $s^\dagger_1 = s^\dagger_2$.
\epf
%-----------------------------------------------------------------------[end proposition]---------------
%-------------------------------------------------------------------------[begin remark]---------------
\begin{rema}
\leavevmode\\
 \ \ \ In the case $l \ne 2$, $-1$ coincides with the unit element $1$ in $(K^\times )^{\wedge}$. 
 Therefore, in the statement of Proposition 2.5, by taking the choice $(C_\delta^\dagger)$ to be such that 1-cocycle map $\delta$ is trivial, we may obtain an ``$s^\dagger$'' satisfying  $s^\dagger \circ \sigma^\dagger = \sigma^\dagger$.
 \end{rema}
%-------------------------------------------------------------------------------[end remark]-------------

%\ere
\leavevmode\\
%%%%%%%%%%%%%%%%%%%%%%%%%%%%%%%%%%%%%%%%%%%%%%%%%
%%%%%%%%%%%%%%%%%%%%%%%%%%%%%---[begin section]---%%%%%%%%%%
\section{The proof of Theorem A}
\leavevmode\\ 
 \ \ \ This section is devoted to proving Theorem A.
 We begin with a review of the notation and setup. Let $l$ be a prime number, $K$ a finite field in which $l$ is invertible, and $\overline{K}$ a separable closure of $K$. We shall denote by $G_K$ the Galois group of $\overline{K}$ over $K$. Next, let $X$ be a hyperbolic curve over $K$ of type ($g_X,r_X$) and
$x^{\mr{log}}$ a strict $K$-rational log point of $\overline{X}^{\mr{log}} := \overline{X}_1^{\mr{log}}$; 
write $\overline{X}_x^{\mr{log}} :=  x^{\mr{log}} \times_{\overline{X}^{\mr{log}}, p_1^\mr{log}} \overline{X}_2^{\mr{log}}$,
$\tilde{x}^{\mr{log}} := x^{\mr{log}} \times_{\overline{X}^{\mr{log}}}  \mbD^{\mr{log}}$. In addition, we assume that we have fixed choices $(C_{i_1}^X)$, $(C_{p_1}^X)$, $(C_{\tilde{x}}^X)$, $(C_\sigma^X)$, $(C_\delta^X)$ [i.e., in the sense that they are \textit{not} subject to conjugacy indeterminacy].

As a first step, we define two actions of $G_K$  on various topological groups, graded Lie algebras, and linear algebraic groups associated to the fundamental groups of $\overline{X}_x^{\mr{log}}$ and $X_2$.  As we shall discuss in the following, these two actions are mapped to one another via the morphisms induced by the switching morphism obtained in Section 2.

%------------------------------------------------------------------------[begin definition]----------------
%\vspace{3mm}
\bde \leavevmode\\
\be
\arrangeVa\item[(i)]
\vspace{-7mm}
The choice $(C_\sigma)$ of a section $\sigma : G_K \migi D_{\tilde{x}}$
 determines,
 by composing with the natural morphism
$D_{\tilde{x}} \migi \Pi_{\overline{X}_x^{\mr{log}}}$  
(resp., $D_{\tilde{x}} \migi \Pi_{X_2}$, $D_{\tilde{x}} \migi \Pi_{X^{\times 2}}$),
a natural action of $G_K$ by conjugation  
on $\Delta_{X_{2/1}}^1 \cong \mr{Ker}(\Pi_{\overline{X}_x^{\mr{log}}} \stackrel{ i_1^\Pi \circ p_1^\Pi}{\migi} \Pi_X)$  (resp.,  $\Delta_{X_2}$,  $\Delta_{X^{\times 2}}$),
hence also on 
\[ \mr{Gr}_{X_{2/1}}^1 := \mr{Gr}_{\mbQ_l}(\Delta_{X_{2/1}}^1 ), \]
\[  \bigl(\text{resp.,} \ \mr{Gr}_{X_2} := \mr{Gr}_{\mbQ_l}(\Delta_{X_2} ), \ \mr{Gr}_{X^{\times 2}} := \mr{Gr}_{\mbQ_l}(\Delta_{X^{\times 2}} ) \bigr), \]
\vspace{-1mm}
\[ \mr{Lie}_{X_{2/1}}^1 := \mr{Lie}(\Delta_{X_{2/1}}^1(1/\infty)), \]
\[ \bigl(\text{resp.,} \ \mr{Lie}_{X_2} := \mr{Lie}(\Delta_{X_2}(1/\infty)), \ \mr{Lie}_{X^{\times 2}} := \mr{Lie}(\Delta_{X^{\times 2}}(1/\infty)) \bigr), \]
\vspace{-1mm}
\[ \mr{Lin}_{X_{2/1}}^1 := \mr{Lin}(\Delta_{X_{2/1}}^1(1/\infty))(\mbQ_l ). \]
\[  \bigl(\text{resp.,} \ \mr{Lin}_{X_2} := \mr{Lin}(\Delta_{X_2}(1/\infty))(\mbQ_l ), \ \mr{Lin}_{X^{\times 2}} := \mr{Lin}(\Delta_{X^{\times 2}}(1/\infty))(\mbQ_l ) \bigr). \]
%\vspace{-1mm}
\\
In the following, we regard these objects as being equipped with these $G_K$-actions.
From the discussion in Definition 1.5 (ii), we have the following commutative diagram consisting of $G_K$-equivariant morphisms
\bcd
\Delta_{X_{2/1}}^1 @> i_1 >> \Delta_{X_2} @> p >> \Delta_{X^{\times 2}}
\\
@VVV @VVV @VVV
\\
\mr{Lin}_{X_{2/1}}^1 @> i_1^{\mr{Lin}} >> \mr{Lin}_{X_2} @> p^{\mr{Lin}} >> \mr{Lin}_{X^{\times 2}}
\ecd
and  topological groups equipped with $G_K$-actions
\[ \Delta^{\mr{Lie}}_{X_2}:= \Delta_{X^{\times 2}} \times_{\mr{Lin}_{X^{\times 2}}} \mr{Lin}_{X_2}, \hspace{5mm} \Pi^{\mr{Lie}}_{X_2} :=  \Delta^{\mr{Lie}}_{X_2} \rtimes G_K \]
as well as   $G_K$-equivariant homomorphisms of topological groups
\[ \mr{Int}^{\Delta}_{X}:\Delta_{X_2} \migi \Delta^{\mr{Lie}}_{X_2}, \hspace{5mm} \mr{Int}^{\Pi}_{X}:\Pi_{X_2} \migi \Pi^{\mr{Lie}}_{X_2}. \]

\arrangeVa\item[(ii)]
Next, the choice $(C_\sigma)$, $(C_\delta)$ yields a new section of the surjective homomorphism $D_{\tilde{x}} \migisurj G_K$ 
\[ \sigma_{\delta} : G_K \longmigi D_{\tilde{x}} \]
\[\hspace{25mm} g \mapsto \delta(g) \cdot \sigma(g)\]
--- which is a homomorphism of topological groups ---.
Then the section $\sigma_{\delta}$ determines, in a similar way to (i), 
a natural action of $G_K$ by conjugation
on 
\[ \breve{\mr{Gr}}_{X_{2/1}}^1 := \mr{Gr}_{\mbQ_l}(\Delta_{X_{2/1}}^1 ), \]
\[  \bigl(\text{resp.,} \ \breve{\mr{Gr}}_{X_2} := \mr{Gr}_{\mbQ_l}(\Delta_{X_2} ), \ \breve{\mr{Gr}}_{X^{\times 2}} := \mr{Gr}_{\mbQ_l}(\Delta_{X^{\times 2}} ) \bigr), \]
\vspace{-1mm}
\[ \breve{\mr{Lie}}_{X_{2/1}}^1 := \mr{Lie}(\Delta_{X_{2/1}}^1(1/\infty)), \]
\[ \bigl(\text{resp.,} \ \breve{\mr{Lie}}_{X_2} := \mr{Lie}(\Delta_{X_2}(1/\infty)), \ \breve{\mr{Lie}}_{X^{\times 2}} := \mr{Lie}(\Delta_{X^{\times 2}}(1/\infty)) \bigr), \]
\vspace{-1mm}
\[ \breve{\mr{Lin}}_{X_{2/1}}^1 := \mr{Lin}(\Delta_{X_{2/1}}^1(1/\infty))(\mbQ_l ). \]
\[  \bigl(\text{resp.,} \ \breve{\mr{Lin}}_{X_2} := \mr{Lin}(\Delta_{X_2}(1/\infty))(\mbQ_l ), \ \breve{\mr{Lin}}_{X^{\times 2}} := \mr{Lin}(\Delta_{X^{\times 2}}(1/\infty))(\mbQ_l ) \bigr). \]
\\
--- where, in the following, we regard these objects as being equipped with the $G_K$-actions 
just defined --- as well as  topological groups equipped with $G_K$-actions
\[ \breve{\Delta}^{\mr{Lie}}_{X_2}:= \Delta_{X^{\times 2}} \times_{\breve{\mr{Lin}}_{X^{\times 2}}} \breve{\mr{Lin}}_{X_2}, \hspace{5mm} \breve{\Pi}^{\mr{Lie}}_{X_2} :=  \breve{\Delta}^{\mr{Lie}}_{X_2} \rtimes G_K .\]

Next, let us recall that by applying  Proposition 2.5, together with the choices $(C_{i_1})$, $(C_{p_1})$, $(C_{\tilde{x}})$, $(C_\delta)$ and the choice $(C_\sigma^\dagger)$ determined naturally by $(C_\sigma)$, we obtain a choice $(C_s^\dagger)$ of a specific automorphism 
$s^\dagger : \Pi_{X_2}^\dagger \migi \Pi_{X_2}^\dagger$.  Let $s^\Pi : \Pi_{X_2} \isom \Pi_{X_2}$ be an automorphism that induces the outer automorphism determined by the switching morphism $s_X : \overline{X}^{\mr{log}}_2 \migi \overline{X}^{\mr{log}}_2$ and is compatible
with $s^\dagger : \Pi_{X_2}^\dagger \migi \Pi_{X_2}^\dagger$.
Then, by Lemma 3.2 below, we obtain $G_K$-equivariant isomorphisms of topological groups
\[ s^{\Delta^{\mr{Lie}}} : \Delta_{X_2}^{\mr{Lie}} \stackrel{\sim}{\longmigi} \breve{\Delta}^{\mr{Lie}}_{X_2}, \hspace{5mm} s^{\Pi^{\mr{Lie}}} : \Pi_{X_2}^{\mr{Lie}} \stackrel{\sim}{\longmigi} \breve{\Pi}^{\mr{Lie}}_{X_2} \]
induced by  $s^\Pi$
and a (non-$G_K$-equivariant) commutative diagram as follows:
\[ \begin{CD}
\Delta_{X_2} @> s^\Delta >> \Delta_{X_2}
\\
@V \mr{Int}_X^\Delta VV @VV \mr{Int}_X^\Delta V
\\
\Delta_{X_2}^{\mr{Lie}} @> s^{\Delta^{\mr{Lie}}} >> \breve{\Delta}^{\mr{Lie}}_{X_2}
\end{CD} 
\hspace{20mm}
\begin{CD}
\Pi_{X_2} @> s^\Pi >> \Pi_{X_2}
\\
@V \mr{Int}_X^\Pi VV @VV \mr{Int}_X^\Pi V
\\
\Pi_{X_2}^{\mr{Lie}} @> s^{\Pi^{\mr{Lie}}} >> \breve{\Pi}^{\mr{Lie}}_{X_2}.
\end{CD}
\]
\ee
\ede
%---------------------------------------------------------------------[end definition]----------------
%-----------------------------------------------------------------------[begin lemma]---------------------
%\vspace{2mm}
\ble
\leavevmode\\
 \ \ \ The $G_K$-action induced by $\sigma_\delta$ (cf. Definition 3.1 (ii)) on $\Delta_{X_2}$ (hence also on $\breve{\mr{Gr}}_{X_2}$, $\breve{\mr{Lie}}_{X_2}$, $\breve{\mr{Lin}}_{X_2}$ and $\breve{\Delta}_{X_2}^{\mr{Lie}}$) coincides with 
the action
 \[ G_K \longmigi \mr{Aut}(\Delta_{X_2}) \]
 \[ \hspace{20mm} g \mapsto   \mr{Inn}(s^{\Pi} \circ i_1^\Pi \circ  \sigma (g)).  \]
 \ele
%\vspace{1mm}
%----------------------------------------------------------------------------[begin proof]---------------
\bpf 
This  follows immediately from condition (3) of Proposition 2.5, together with the definition of the $G_K$-action induced by $\sigma_\delta$.\epf
%---------------------------------------------------------------------------[end lemma]------------------
%-----------------------------------------------------------------------[begin lemma]---------------------
%\vspace{2mm}
\ble
\leavevmode\\
 \ \ \ $\mr{Int}^{\Delta}_{X}$ and $\mr{Int}^{\Pi}_{X}$ are injective
 (cf. ~\cite{HSHcusp}, Lemma 4.3 in the case where $X$ is proper).
\ele
%\vspace{1mm}
%----------------------------------------------------------------------------[begin proof]---------------
\bpf 
It suffices to verify that  $\Delta_{X_2} \migi \mr{Lin}_{X_2}$ is injective.
But this follows from the discussion in Definition 1.5 (ii) and the fact that
$ \bigcap_{m \ge 1} \Delta_{X^{\times 2}} (m) = 1 $
(cf. ~\cite{STAB}, Corollary 2.6).   
\epf
%---------------------------------------------------------------------------[end lemma]------------------
%\vspace{3mm}

Next, we shall construct certain graded Lie algebras equipped with a $G_K$-action --- which we shall denote by $\mcL^1_X$ and $\mcL^2_X$ --- by using various subgroups of $\Pi_{\overline{X}^{\mr{log}}_x}$. Comparing these graded Lie algebras to the graded Lie algebras discussed above (cf. Lemma 3.5, 3.6) will allow us to reconstruct various groups  associated to $\Pi_{X_2}$ from those associated to $\Pi_{\overline{X}^{\mr{log}}_x}$ (cf. Proposition 3.8).
This will play an important role in the proof of Theorem A.
%-------------------------------------------------------------------------[begin definition]------------
%\vspace{2mm}
\bde\leavevmode\\
\be
\arrangeVa\item[(i)]
\vspace{-7mm}
For each $j=1,2,\dotsm,r$, let us fix a choice of the inertia subgroup $I_j \subseteq \Delta_{X_{2/1}}^1 \cong \mr{Ker}(\Pi_{\overline{X}_x^{\mr{log}}} \stackrel{p^\Pi_1 \circ i^\Pi_1}{\migi} \Pi_X)$ associated to the $j$-th cusp (relative to some ordering of the cusps of $X \times_K \overline{K}$) among the various $\Delta_{X_{2/1}}^1$-conjugates of these subgroups.
Then, we have canonical isomorphisms 
\[ \hspace{10mm} \eta_j : I_{\tilde{x}} \stackrel{\sim}{\longmigi} I_j  \ \  (j= 1,2, \dotsm , r)\]
Indeed,  recall that the {\it kernel} of the natural quotient $(\Delta_{X_{2/1}}^1)^{\mr{ab}} \migisurj \Delta_{\overline{X}}^{\mr{ab}}$ coincides with the submodule
$\bigoplus_{j'=1}^r I_{j'} \subseteq (\Delta_{X_{2/1}}^1)^{\mr{ab}}$;
thus, since the subgroup $I_{\tilde{x}}$ of $(\Delta_{X_{2/1}}^1)^{\mr{ab}}$
is contained in this kernel, it follows that the composite
$I_{\tilde{x}} \migiincl \bigoplus^r_{j'=1} I_{j'} \migisurj I_j \stackrel{(-1)}{\migi} I_j$
of this inclusion with the natural projection to $j$-th factor multiplied by $-1$ 
 yields the required  isomorphism.

 For $n=1,2$ we shall denote by $\mcV^n$ the completion with respect to the filtration topology of the free Lie algebra generated by
\[ V^n :=  I_{\tilde{x}}  \oplus \bigl( \bigoplus_{j =1}^r  I_j  \oplus \Delta_{\overline{X}}^{\mr{ab}} \bigr) ^{\oplus n}  \] 
equipped with a natural grading (hence also a filtration) by taking $I_{\tilde{x}}, I_j$ to be of weight 2, $\Delta_{\overline{X}}^{\mr{ab}}$ to be of weight 1. 

\arrangeVa\item[(ii)]
If $X$ has \textit{genus} $\geq 1$, then we shall write
\[ M_X := \mr{Hom}_{\mbZ_l}(H^2(\Delta_{\overline{X}},\mbZ_l), \mbZ_l). \]
Note that $M_X$ is canonically isomorphic to $I_{\tilde{x}}$ as a $G_K$-module.
Indeed, recall the natural quotient
$(\Delta_{X_{2/1}}^1 /<I_j>_{j=1, \dotsm , r})\migisurj \Delta_{\overline{X}}$;
the associated maximal cuspidally central quotient (cf.  ~\cite{MZKcusp}, Definition 1.1 (i))
yields an extension of $\Delta_{\overline{X}}$ by $I_{\tilde{x}}$; this extension determines a generator of the rank one free $\mbZ_l$-module $H^2(\Delta_{\overline{X}}, I_{\tilde{x}}) \cong \mr{Hom}_{\mbZ_l}(M_X, I_{\tilde{x}})$ (cf., e.g.,  ~\cite{MZKsect}, Lemma 4.2, (i), (ii), (iii) ), hence an isomorphism $M_X \isom  I_{\tilde{x}}$, as desired.

The cup product on the group cohomology of $\Delta_{\overline{X}}$
\[ \bigwedge^2 H^1(\Delta_{\overline{X}}, M_X) \longmigi H^2(\Delta_{\overline{X}}, M_X \otimes_{\mbZ_l} M_X) \cong M_X \]
determines an isomophism
\[ \bigl( H^1(\Delta_{\overline{X}}, M_X) \cong \bigr)\  \mr{Hom}(\Delta_{\overline{X}}^{\mr{ab}}, M_X) \ \stackrel{\sim}{\longmigi}\  \Delta_{\overline{X}}^{\mr{ab}}\  \bigl( \cong \mr{Hom}(H^1(\Delta_{\overline{X}}, M_X), M_X) \bigr), \]
hence  composites of natural homomorphisms
\[ \phi : I_{\tilde{x}} \stackrel{\sim}{\longmigi} M_X \longmigi \bigwedge^2 {\Delta_{\overline{X}}^{\mr{ab}}}, \hspace{5mm}  \psi : \bigwedge^2 {\Delta_{\overline{X}}^{\mr{ab}}} \longmigi M_X \stackrel{\sim}{\longmigi}  I_{\tilde{x}}. \]

 If $X$ has \textit{genus} $0$, then we take $\phi$, $\psi$ to be the zero maps.

%----------------------------------------------------------------------------------------------------------
\arrangeVa\item[(iii)]
We define $\mcL_X^n $ to be the quotient of $\mcV^n$ by the relations determined by the images of the following morphisms (which are patterned after the presentations given in Proposition 1.9):
\\[0mm]
\begin{itemize}
%----------------------------------------------------------------------------------------------------------
\item[(1)] 
When $n=1$, \vspace{2mm}
\begin{itemize} 
\item[$\bullet_1$] $I_{\tilde{x}} \longmigi \mcV^1 (2/3) \ ; \ m \mapsto (\text{id}_{I_{\tilde{x}}} + \sum \eta_j + \phi )(m)$ \vspace{3mm}
\end{itemize}
%---------------------------------------------------------------------------------------------------------
%\vspace{0mm}

\item[(2)] When $n=2$ $(1 \leq i \leq g, \  1 \leq j, j'  \leq r, \ j \neq j', \ \{k , k' \} = \{1, 2 \} )$,
\begin{itemize} \vspace{2mm}
\item[$\bullet_1$] $I_{\tilde{x}} \longmigi \mcV^2 (2/3) \ ; \ m \mapsto m+ i_k (\sum \eta_j + \phi )(m)$ \vspace{3mm}%\notag \\
\item[$\bullet_2$]  
$I_{\tilde{x}} \otimes_{\mbZ_l} \Delta_{\overline{X}}^{\mr{ab}} \longmigi \mcV^2 (3/4) \ ; m \otimes a \mapsto  [ i_k \circ \eta_j (m),  i_{k'} (a)]$ \vspace{3mm}
\item[$\bullet_3$] $I_{\tilde{x}} \longmigi \mcV^2 (4/5) \ ; \ m \mapsto [i_k \circ \eta_j (m), i_{k'} \circ \eta_{j'}(m)]$ \vspace{3mm}  

\item[$\bullet_{4,5}$] $\bigwedge^2 {\Delta_{\overline{X}}^{\mr{ab}}} \longmigi \mcV^2 (2/3) \ ; a \wedge a'  \mapsto [ i_k(a) , i_{k'} (a') ] -\psi (a \wedge a')$ \vspace{3mm}%\notag \\
\end{itemize}
\end{itemize}
--- where $``[ \ , \ ]"$ denotes the Lie bracket, and for $k=1,2$, \,$i_k : ( \bigoplus I_j \oplus \bigwedge^2 {\Delta_{\overline{X}}^{\mr{ab}}}) \migiincl ( \bigoplus I_j \oplus \bigwedge^2 {\Delta_{\overline{X}}^{\mr{ab}}})^{\oplus 2} $ denotes the inclusion into the $k$-th factor.
\\[-1mm]

\arrangeVa\item[(iv)]
The natural $G_K$-action on each direct summand in $\mcV^n$ determines a natural
 $G_K$-action on $\mcV^n$.  One verifies immediately that the ideal generated by the
 relations defined in (iii) is {\it preserved} by this $G_K$-action.
Thus, we obtain a natural $G_K$-action on the graded Lie algebra 
\[ \mcL_X^1   \ \  (\text{resp.}, \mcL_X^2) \]
and a $G_K$-equivariant homomorphism 
\[ i_1^\mcL : \mcL_X^1 \longmigi \mcL_X^2 \]
 of graded Lie algebras determined by the map on generators given by
\[ I_{\tilde{x}}  \oplus \big( \bigoplus_{j =1}^r  I_j  \oplus \Delta_{\overline{X}}^{\mr{ab}} \big) \longmigi I_{\tilde{x}}  \oplus \bigl( \bigoplus_{j =1}^r  I_j  \oplus \Delta_{\overline{X}}^{\mr{ab}} \bigr) ^{\oplus 2} \]
\\[-8mm]
\[ \ ( a, b ) \hspace{5mm} \mapsto \hspace{5mm} ( a, i_1(b) ) ,\  \]
as well as a $G_K$-equivariant  isomorphism
\[ s_X^\mcL : \mcL_X^2 \stackrel{\sim}{\longmigi} {\mcL}_X^2 \]
of graded Lie algebras determined by the map on generators given by
\[ I_{\tilde{x}}  \oplus \bigl( \bigoplus_{j =1}^r  I_j  \oplus \Delta_{\overline{X}}^{\mr{ab}} \bigr) ^{\oplus 2} \longmigi I_{\tilde{x}}  \oplus \bigl( \bigoplus_{j =1}^r  I_j  \oplus \Delta_{\overline{X}}^{\mr{ab}} \bigr) ^{\oplus 2} \]
\[ (a, b_1, b_2) \hspace{5mm} \mapsto  \hspace{5mm} (a, b_2, b_1) .\]
\ee
\ede
%----------------------------------------------------------------------------[end definition]-------------
%---------------------------------------------------------------------------[begin lemma]---------
%\vspace{3mm}
\ble\leavevmode\\
 \ \ \ Consider the homomorphism of graded Lie algebras
 $\mcV^1 \migi \mr{Gr}_{\mbQ_l}(\Delta^1_{X_{2/1}})$
determined by the natural inclusions $\Delta_{\overline{X}}^{\mr{ab}} \migiincl \mr{Gr}_{\mbQ_l}(\Delta^1_{X_{2/1}})(1/2)$, $I_{\tilde{x}} \migiincl \mr{Gr}_{\mbQ_l}(\Delta^1_{X_{2/1}})(2/3)$ and $I_j \migiincl \mr{Gr}_{\mbQ_l}(\Delta^1_{X_{2/1}})(2/3)$.
This homomorphism of graded Lie algebras factors through $\mcL_X^1$, and the resulting homomorphism
 $\mfh^1 : \mcL_X^1 \migi \mr{Gr}_{\mbQ_l}(\Delta^1_{X_{2/1}})$ is a $G_K$-equivariant isomorphism of graded Lie algebras, whether we regard $\mr{Gr}_{\mbQ_l}(\Delta^1_{X_{2/1}})$ as the underlying graded Lie algebra (i.e., without $G_K$-action) of $\mr{Gr}_{X_{2/1}}^1$ or as the underlying graded Lie algebra of $\breve{\mr{Gr}}_{X_{2/1}}^1$.
\ele
%-----------------------------------------------------------------------------[begin proof]---------------
%\vspace{0mm}
\bpf 
The asserted {\it $G_K$-equivariance} follows immediately from the definitions.  Thus, it suffices to verify that $\mfh^1$ is an {\it isomorphism}.
When $x$ is {\it not} a cusp of $X$, this follows immediately from Proposition 1.9 (i), applied to $X_x$.  Thus, it suffices to verify that $\mfh^1$ is an isomorphism in the case where $x$ is a cusp of $X$.
Let  $S$ be a mixed characteristic trait (i.e., the spectrum of a henselian discrete valuation ring) whose residue field is isomorphic to the residue field of $x$, and write $S^{\mr{log}}$ for the log scheme obtained by equipping $S$ with the log structure determined by the closed point of $S$. Next, let us take a stable log curve 
$\overline{X}_S^{\mr{log}} \migi S^{\mr{log}}$ 
 whose special fiber is isomorphic to $\overline{X}_x^{\mr{log}} \migi x^{\mr{log}}$
 and such that the interior $U$ of $X_S^{\mr{log}}$ is a hyperbolic curve over the fraction field of $S$.
Then (cf. the discussion of ~\cite{MZKcomb2}, \S 0, in the characteristic zero case) we obtain a  natural isomorphism
$\Delta_{\overline{X}_x^{\mr{log}}} \isom \Delta_U$
by composing
a certain ``\textit{specialization isomorphism}"
$\Delta_{\overline{X}_x^{\mr{log}}} \isom \Delta_{\overline{X}_S^{\mr{log}}}$
with an isomorphism $\Delta_{\overline{X}_S^{\mr{log}}} \isom \Delta_U$
arising from the ``\textit{log purity theorem}".
Hence, the fact that $\mfh^1$ is an isomorphism follows immediately by applying this isomorphism 
$\Delta_{\overline{X}_x^{\mr{log}}} \isom \Delta_U$, together with Proposition 1.9 (i). 
\epf
%----------------------------------------------------------------------------[end lemma]-------------

%---------------------------------------------------------------------------[begin lemma]---------
%\vspace{3mm}
\ble\leavevmode\\
 \ \ \ Let 
\[ i_1^{\mr{Lie}} : \mr{Lie}_{X_{2/1}}^1 \longmigi \mr{Lie}_{X_2}, \ \ \breve{i}_1^{\mr{Lie}} : \breve{\mr{Lie}}_{X_{2/1}}^1 \longmigi \breve{\mr{Lie}}_{X_2}, \ \   s_{X}^{\mr{Lie}} : \mr{Lie}_{X_{2}} \stackrel{\sim}{\longmigi} \breve{\mr{Lie}}_{X_2} \]
 be the $G_K$-equivariant homomorphisms of graded Lie algebras induced by
 $i_1^\Pi : \Pi_{\overline{X}^{\mr{log}}_x} \migi \Pi_{X_2}$, $i_1^\Pi : \Pi_{\overline{X}^{\mr{log}}_x} \migi \Pi_{X_2}$
 and $s_{X}^\Pi : \Pi_{X_2} \isom \Pi_{X_2}$ respectively. 

Then there exist $G_K$-equivariant isomorphisms of graded Lie algebras
\[ h_X^1 : \mcL^1_X \stackrel{\sim}{\longmigi} \mr{Lie}_{X_{2/1}}^1, \ \ \ \breve{h}_X^1 : \mcL^1_X \stackrel{\sim}{\longmigi} \breve{\mr{Lie}}_{X_{2/1}}^1, \]
\[ h_X^2 : \mcL^2_X \stackrel{\sim}{\longmigi} \mr{Lie}_{X_2}, \ \ \ \breve{h}_X^2 : {\mcL}^2_X \stackrel{\sim}{\longmigi} \breve{\mr{Lie}}_{X_2} \]
which fit into the following commutative diagrams consisting of $G_K$-equivariant mophisms
\[ \begin{CD}
\mcL^1_X @> i_1^\mcL >> \mcL^2_X 
\\
@V h_X^1 V \wr V @V h_X^2 V \wr V
\\
\mr{Lie}_{X_{2/1}}^1 @> i_1^{\mr{Lie}} >> \mr{Lie}_{X_2}
\end{CD} \hspace{10mm}
\begin{CD}
\mcL^1_X @> i_1^\mcL >> \mcL^2_X 
\\
@V \breve{h}_X^1 V \wr V @V \breve{h}_X^2 V \wr V
\\
\breve{\mr{Lie}}_{X_{2/1}}^1 @> \breve{i}_1^{\mr{Lie}} >> \breve{\mr{Lie}}_{X_2}
\end{CD} \hspace{10mm}
\begin{CD}
\mcL^2_X @> s_{X}^\mcL >> {\mcL}^2_X 
\\
@V h_X^2 V \wr V @V \breve{h}_X^2 V \wr V
\\
\mr{Lie}_{X_{2}} @> s_{X}^{\mr{Lie}} >> \breve{\mr{Lie}}_{X_2}.
\end{CD} \]
\ele
%-----------------------------------------------------------------------------[begin proof]---------------
%\vspace{3mm}
\bpf 
Consider the homomorphism of graded Lie algebras
 $\mcV^2 \migi \mr{Gr}_{\mbQ_l}(\Delta_{X_2})$
determined by
\[ I_{\tilde{x}}  \oplus \bigl( \bigoplus_{j =1}^r  I_j  \oplus \Delta_{\overline{X}}^{\mr{ab}} \bigr) ^{\oplus 2} \longmigi \mr{Gr}_{\mbQ_l}(\Delta_{X_2}) \hspace{45mm} \]
\[ \hspace{0mm} (a, b_1, b_2) \hspace{5mm} \mapsto  \ \ i^\Pi_1 (a +b_1) + s^\Pi_{X} \circ i_1^{\Pi} (b_2). \]
Then it follows from Proposition 1.9 (ii) that this homomorphism of graded Lie algebras factors through $\mcL_X^2$, and that the resulting homomorphism
$\mfh^2 : \mcL_X^2 \migi \mr{Gr}_{\mbQ_l}(\Delta_{X_2})$
 is a $G_K$-equivariant isomorphism of graded Lie algebras, whether we regard
$\mfh^2$
as the morphism of underlying graded Lie algebras (i.e., without $G_K$-actions) $h_{}^2 : \mcL^2_X \migi \mr{Gr}_{X_2}$
or as the morphism of underlying graded Lie algebras 
$\breve{h}_{}^2 : {\mcL}^2_X \migi \breve{\mr{Gr}}_{X_2}$.
If we denote by
$ i_1^{\mr{Gr}} : \mr{Gr}_{X_{2/1}}^1 \migi \mr{Gr}_{X_2}$,
$ \breve{i}_1^{\mr{Gr}} : \breve{\mr{Gr}}_{X_{2/1}}^1 \migi \breve{\mr{Gr}}_{X_2}$,
$s_{X}^{\mr{Gr}} : \mr{Gr}_{X_{2}} \stackrel{\sim}{\migi} \breve{\mr{Gr}}_{X_2}$
 the $G_K$-equivariant homomorphisms of graded Lie algebras induced by
$i_1^\Pi : \Pi_{\overline{X}^{\mr{log}}_x} \migi \Pi_{X_2}$,
$i_1^\Pi : \Pi_{\overline{X}^{\mr{log}}_x} \migi \Pi_{X_2}$,
and
$s_{X}^\Pi : \Pi_{X_2} \isom \Pi_{X_2}$,
respectively,
then we obtain $G_K$-equivariant commutative diagrams as follows:
\[ \begin{CD}
\mcL^1_X @> i_1^\mcL >> \mcL^2_X 
\\
@V h_{}^1 V \wr V @V h_{}^2 V \wr V
\\
\mr{Gr}_{X_{2/1}}^1 @> i_1^{\mr{Gr}} >> \mr{Gr}_{X_2}
\end{CD} \hspace{10mm}
\begin{CD}
\mcL^1_X @> i_1^\mcL >> \mcL^2_X 
\\
@V \breve{h}_{}^1 V \wr V @V \breve{h}_{}^2 V \wr V
\\
\breve{\mr{Gr}}_{X_{2/1}}^1 @> \breve{i}_1^{\mr{Gr}} >> \breve{\mr{Gr}}_{X_2}
\end{CD} \hspace{10mm}
\begin{CD}
\mcL^2_X @> s_X^\mcL >> {\mcL}^2_X 
\\
@V h_{}^2 V \wr V @V \breve{h}_{}^2 V \wr V
\\
\mr{Gr}_{X_{2}} @> s_X^{\mr{Gr}} >> \breve{\mr{Gr}}_{X_2}.
\end{CD} \]
On the other hand, it follows from Proposition 1.8 that  we have $G_K$-equivariant commutative diagrams as follows:
\[ \begin{CD}
\mr{Gr}_{X_{2/1}}^1 @> i_1^{\mr{Gr}} >> \mr{Gr}_{X_2}
\\
@V  V \wr V @V V \wr V
\\
\mr{Lie}_{X_{2/1}}^1 @> i_1^{\mr{Lie}} >> \mr{Lie}_{X_2}
\end{CD} \hspace{10mm}
\begin{CD}
\breve{\mr{Gr}}_{X_{2/1}}^1 @> \breve{i}_1^{\mr{Gr}} >> \breve{\mr{Gr}}_{X_2}
\\
@V  V \wr V @V V \wr V
\\
\breve{\mr{Lie}}_{X_{2/1}}^1 @> \breve{i}_1^{\mr{Lie}} >> \breve{\mr{Lie}}_{X_2}
\end{CD} \hspace{10mm}
\begin{CD}
\mr{Gr}_{X_{2}} @> s_X^{\mr{Gr}} >> \breve{\mr{Gr}}_{X_2}.
\\
@V  V \wr V @V  V \wr V
\\
\mr{Lie}_{X_{2}} @> s_X^{\mr{Lie}} >> \breve{\mr{Lie}}_{X_2}.
\end{CD} \]
By composing the vertical arrows in these commutative diagrams, we obtain the required isomorphisms.
\epf
%----------------------------------------------------------------------------[end lemma]-------------
%\vspace{3mm}

Now, let $L$ be a finite field of cardinality prime to $l$,
$Y$ a hyperbolic curve over $L$ of type ($g_Y, r_Y$),
$y^{\mr{log}}$ a strict $L$-rational log point of
$\overline{Y}^{\mr{log}} := \overline{Y}_1^{\mr{log}}$;
we shall use similar notation for objects obtained from $Y$
(e.g., $Y_2$, $\overline{Y}_y^{\mr{log}}$, $\Pi_{Y_2}$, $\Pi_{\overline{Y}_y^{\mr{log}}}$, etc.) to the notation used for objects obtained from $X$. 

%\vspace{3mm}
%---------------------------------------------------------------------------[begin definition]-----------
\bde\leavevmode\\
 \begin{itemize} \vspace{-5mm} \item[(i)]
Consider an isomorphism of profinite groups $\alpha : \Pi_X \isom \Pi_Y$ (resp., $\Pi_{\overline{X}_x^{\mr{log}}} \isom \Pi_{\overline{Y}_y^{\mr{log}}}$).
Then the natural surjections 
$\Pi_X \migisurj G_K$,
$ \Pi_Y \migisurj G_L$
 (resp., $\Pi_{\overline{X}_x^{\mr{log}}} \migisurj G_{K} $,
 $\Pi_{\overline{Y}_x^{\mr{log}}} \migisurj G_{L}$
 ) arising from the structure morphisms over  finite fields  
may be characterized group-theoretically (cf. ~\cite{TAM}, Proposition 3.3) as the (unique) maximal ($\hat{\mbZ}$-)free abelian quotient.
 Thus, $\alpha$ induces an isomorphism $G_K \isom G_L$.

We shall say that $\alpha$ is \textit{Frobenius-preserving} if the isomorphism $G_K \isom G_L$ obtained as above preserves the Frobenius elements. 

\item[(ii)]
We shall denote by 
\[ (C^{X \isom Y}_{x,y}) \ (\text{resp.,} \ (C^{\overline{X}^{\mr{log}}_x \isom \overline{Y}^{\mr{log}}_y}_{x,y}) ) \]
 a choice of a specific Frobenius-preserving isomorphism $\Pi_X \isom \Pi_Y$ (resp., $\Pi_{\overline{X}_x^{\mr{log}}} \isom \Pi_{\overline{Y}_y^{\mr{log}}}$)
which maps the decomposition group of $x$ (resp., the diagonal cusp $\tilde{x}$) onto the decomposition group of $y$ (resp., the diagonal cusp $\tilde{y}$) up to conjugation.
\end{itemize}
\ede
%-------------------------------------------------------------[begin proposition]------------------
%\vspace{3mm}
\bpr\leavevmode\\
 \ \ \ Let us fix specific choices of 
$(C^Y_{i_1})$,  $(C^Y_{p_1})$, $(C^Y_{\tilde{y}})$ and $(C^{\overline{X}^{\mr{log}}_x \isom \overline{Y}^{\mr{log}}_y}_{x,y})$.
Denote by $\alpha:\Pi_{\overline{X}_x^{\mr{log}}} \isom \Pi_{\overline{Y}_y^{\mr{log}}}$
the isomorphism determined by $(C^{\overline{X}^{\mr{log}}_x \isom \overline{Y}^{\mr{log}}_y}_{x,y})$.
Let us assume that  the decomposition subgroups determined by  $(C^X_{\tilde{x}})$, $(C^Y_{\tilde{y}})$ are compatible with respect to $\alpha$. 
\begin{itemize} \vspace{0mm}
\item[(i)]
There exists a unique pair consisting of a choice $(C^Y_\theta)$
of a specific 1-cocycle
$\theta : G_L \migi I_{\tilde{y}} :=  \mr{Ker}( D_{\tilde{y}} \migisurj \Pi_Y)$
and a choice 
 $(C^Y_\tau)$
of a specific section $\tau : G_L \migi D_{\tilde{y}}$
which are compatible with
$(C^X_\delta)$ and $(C^X_\sigma)$, respectively, via $(C^{X \isom Y}_{x,y})$
in an evident fashion.
 \item[(ii)]There exists a $G_K$-equivariant  isomorphism
$\alpha^{\Delta^{\mr{Lie}}}_2:\Delta^{\mr{Lie}}_{X_2} \isom \Delta^{\mr{Lie}}_{Y_2}$
of topological groups satisfying the following conditions:
\begin{itemize}
%----------------------------------------------------------------------------------------------------------
\item[(1)] 
$\alpha^{\Delta^{\mr{Lie}}}_2$  is also $G_K$-equivariant when we regard it as a map $\breve{\Delta}^{\mr{Lie}}_{X_2} \isom \breve{\Delta}^{\mr{Lie}}_{Y_2}$ under the natural identifications $\Delta^{\mr{Lie}}_{X_2} \isom \breve{\Delta}^{\mr{Lie}}_{X_2}$, $\Delta^{\mr{Lie}}_{Y_2} \isom \breve{\Delta}^{\mr{Lie}}_{Y_2}$ without $G_K$-actions.
\item[(2)]
If we denote by $\alpha^{\Pi^{\mr{Lie}}}_2:\Pi^{\mr{Lie}}_{X_2} \isom \Pi^{\mr{Lie}}_{Y_2}$, 
$\breve{\alpha}^{\Pi^{\mr{Lie}}}_2:\breve{\Pi}^{\mr{Lie}}_{X_2} \isom \breve{\Pi}^{\mr{Lie}}_{Y_2}$
the semi-direct products of $\alpha^{\Delta^{\mr{Lie}}}_2$ with the isomorphism
$G_K\isom G_L$ (i.e., determined by $\alpha$) relative to the respective actions of $G_K$ and $G_L$, then these morphisms
 make the following diagrams commute
\[ \begin{CD}
\Pi_{\overline{X}_x^{\mr{log}}} @> \mr{Int}_X^\Pi \circ i_1^\Pi >>\Pi^{\mr{Lie}}_{X_2} @>>> \Pi_{X^{\times 2}}
\\
@V \alpha V \wr V @V \alpha^{\Pi^{\mr{Lie}}}_2 V \wr V @V \overline{\alpha} \times \overline{\alpha} V \wr V
\\
\Pi_{\overline{Y}_y^{\mr{log}}} @> \mr{Int}_Y^\Pi \circ i_1^\Pi >>\Pi^{\mr{Lie}}_{Y_2} @>>> \Pi_{Y^{\times 2}}
\end{CD}
 \ \ \ \ \ \ \ \ \ \ 
\begin{CD}
\Pi^{\mr{Lie}}_{X_2} @> s_{X}^{\Pi^{\mr{Lie}}} >> \breve{\Pi}^{\mr{Lie}}_{X_2}
\\
@V \alpha^{\Pi^{\mr{Lie}}}_2 VV @V \breve{\alpha}^{\Pi^{\mr{Lie}}}_2 VV
\\
\Pi^{\mr{Lie}}_{Y_2} @> s_{Y}^{\Pi^{\mr{Lie}}} >> \breve{\Pi}^{\mr{Lie}}_{Y_2}.
\end{CD} \]
\end{itemize}
\end{itemize}
\epr
%------------------------------------------------------------[begin proof]------------------------
%\vspace{2mm}
\bpf
Assertion (i) follows immediately by ``transport of structure".
Next, we consider assertion (ii).
Since $\alpha$ is assumed to be Frobenius-preserving, it follows from  ~\cite{MZKcomb1}, Corollary 2.7 (i) that $(g_X,r_X) = (g_Y,r_Y)$, and that $\alpha$ induces an isomorphism $\alpha^{\mr{cpt}} : \Delta^{\mr{ab}}_{\overline{X}} \isom \Delta^{\mr{ab}}_{\overline{Y}}$ and a bijective correspondence between the respective sets of cusps of $\overline{X}^{\mr{log}}_x, \overline{Y}_y^{\mr{log}}$ as well as isomorphisms of the inertia subgroups of cusps corresponding via this bijection. 
By applying these isomorphisms (together with the construction of $\mcL_X^1$, $\mcL_X^2$, $\mcL_Y^1$, $\mcL_Y^2$), Lemma 3.6 yields $G_K$-equivariant isomorphisms
$\alpha^{\mr{Lie}} : \mr{Lie}_{X_{2/1}}^1 \cong \mr{Lie}_{Y_{2/1}}^1$,
$\breve{\alpha}^{\mr{Lie}} : \breve{\mr{Lie}}_{X_{2/1}}^1 \cong \breve{\mr{Lie}}_{Y_{2/1}}^1$, $\alpha_2^{\mr{Lie}} : \mr{Lie}_{X_{2}} \cong \mr{Lie}_{Y_{2}}$ and $\breve{\alpha}_2^{\mr{Lie}} : \breve{\mr{Lie}}_{X_{2}} \cong \breve{\mr{Lie}}_{Y_{2}}$.
These morphisms give rise to a $G_K$-equivariant commutative diagram as follows:
\bcd
\mr{Lie}_{X_{2/1}}^1 @> i_1^{\mr{Lie}} >> \mr{Lie}_{X_2} @> s_{X}^{\mr{Lie}} >> \breve{\mr{Lie}}_{X_2} @<  \breve{i}_1^{\mr{Lie}} << \breve{\mr{Lie}}_{X_{2/1}}^1
\\
@V \alpha^{\mr{Lie}} V \wr V @V \alpha_2^{\mr{Lie}} V\wr V @V \breve{\alpha}_2^{\mr{Lie}} V \wr V @V \breve{\alpha}^{\mr{lie}} V \wr V 
\\
\mr{Lie}_{Y_{2/1}}^1 @> i_1^{\mr{Lie}} >> \mr{Lie}_{Y_2} @> s_{Y}^{\mr{Lie}} >> \breve{\mr{Lie}}_{Y_2} @< \breve{i}_1^{\mr{Lie}} <<  \breve{\mr{Lie}}_{Y_{2/1}}^1.
\ecd
Then it follows from the {\it functoriality} of $\mr{Lin}( - )$
that  we obtain a $G_K$-equivariant commutative diagram as follows:  
\bcd
\mr{Lin}_{X_{2/1}}^1 @> i_1^{\mr{Lin}} >> \mr{Lin}_{X_2} @> s_{X}^{\mr{Lin}} >> \breve{\mr{Lin}}_{X_2} @<  \breve{i}_1^{\mr{Lin}} << \breve{\mr{Lin}}_{X_{2/1}}^1
\\
@V \alpha^{\mr{Lin}} V \wr V @V \alpha_2^{\mr{Lin}} V\wr V @V \breve{\alpha}_2^{\mr{Lin}} V \wr V @V \breve{\alpha}^{\mr{Lin}} V \wr V 
\\
\mr{Lin}_{Y_{2/1}}^1 @> i_1^{\mr{Lin}} >> \mr{Lin}_{Y_2} @> s_{Y}^{\mr{Lin}} >> \breve{\mr{Lin}}_{Y_2} @< \breve{i}_1^{\mr{Lin}} <<  \breve{\mr{Lin}}_{Y_{2/1}}^1.
\ecd
Note (cf. ~\cite{MZKcusp}, Remark 35)
that modifying the choice $(C_\sigma^X)$ of a specific section
$G_K \migi D_{\tilde{x}}$
by a cocycle $G_K\migi I_{\tilde{x}}$ determined by the choice $(C_\delta^X)$
affects the Galois invariant isomorphisms of Proposition 1.8, (ii),
by conjugation by an element  $c_X$ of the subgroup
obtained by tensoring  $I_{\tilde{x}}$  with $\mbQ_l$;
a similar statement holds, with respect to some ``$c_Y$'',  for objects associated to $Y$
when we modify $(C_\tau^Y)$ by $(C_\theta^Y)$. 
One may verify easily that $\alpha$ maps $c_X$ to $c_Y$,
hence that $\alpha^{\mr{Lin}} =  \breve{\alpha}^{\mr{Lin}}$  as a morphism of underlying topological groups (i.e., without  $G_K$-actions).
Next, recall that the morphisms
$i_1^{\mr{Lin}}$ and $s_{X}^{\mr{Lin}} \circ i_1^{\mr{Lin}}$
are compatible with the corresponding morphisms ``$i_1^{\mr{Lin}}$'' and ``$s_{Y}^{\mr{Lin}} \circ i_1^{\mr{Lin}}$'' associated to $Y$
via the natural identification of $\mr{Lin}_{X_2}$ with $\breve{\mr{Lin}}_{X_2}$
(i.e., without  $G_K$-actions).
Also, let us recall that $\mr{Lin}_{X_2}$
(resp., $\breve{\mr{Lin}}_{X_2}$)
is generated by the {\it images} of 
$\mr{Lin}_{X_{2/1}}^1 \stackrel{i_1^{\mr{Lin}}}{\migi} \mr{Lin}_{X_2}$
(resp., $\breve{\mr{Lin}}_{X_{2/1}}^1 \stackrel{\breve{i}_1^{\mr{Lin}}}{\migi} \breve{\mr{Lin}}_{X_2}$)
and the composite
$\mr{Lin}_{X_{2/1}}^1 \stackrel{i_1^{\mr{Lin}}}{\migi} \mr{Lin}_{X_2} \stackrel{s_{X}^{\mr{Lin}}}{\migi} \breve{\mr{Lin}}_{X_2} = {\mr{Lin}}_{X_2}$
(resp., $\breve{\mr{Lin}}_{X_{2/1}}^1 \stackrel{\breve{i}_1^{\mr{Lin}}}{\migi} \breve{\mr{Lin}}_{X_2} = {\mr{Lin}}_{X_2} \stackrel{s_{X}^{\mr{Lin}}}{\migi} \breve{\mr{Lin}}_{X_2}$).
Since the restrictions of $\alpha_2^{\mr{Lin}}$ and $\breve{\alpha}_2^{\mr{Lin}}$ to these {\it image} subgroups 
coincide by virtue of  the equality $\alpha^{\mr{Lin}} =  \breve{\alpha}^{\mr{Lin}}$,
 we obtain  that $\alpha_2^{\mr{Lin}} =\breve{\alpha}_2^{\mr{Lin}}$.
Therefore, by construction, $\alpha_2^{\mr{Lin}}$($ =\breve{\alpha}_2^{\mr{Lin}}$)
induces the required
$G_K$-equivariant  isomorphism
$\alpha^{\Delta^{\mr{Lie}}}_2:\Delta^{\mr{Lie}}_{X_2} \isom \Delta^{\mr{Lie}}_{Y_2}$
of topological groups satisfying conditions (1), (2).
This completes the proof of assertion (ii).
\epf
%---------------------------------------------------------------[end proposition]----------------

One of main results of this paper, i.e., (a slightly generalized version of) Theorem A, is the following:

%-----------------------------------------------------------------[begin theorem]---------------
%\vspace{2mm}
\bt\leavevmode\\
 \ \ \ Let $X$ (\text{resp.}, $Y$) be a hyperbolic curve over a finite field $K$ (\text{resp.}, $L$),
$x$ a $K$-rational point of $\overline{X}$ (\text{resp.}, $y$ an $L$-rational point of $\overline{Y}$),
$X_2$ (\text{resp.}, $Y_2$) the second configuration space associated to $X$ (\text{resp.}, $Y$), $\overline{X}_x^{\mr{log}}$ (\text{resp.}, $\overline{Y}_y^{\mr{log}}$) the cuspidalization of $X$ at $x$ (\text{resp.}, of $Y$ at $y$) [cf. Definition 2.2],
$D_{\tilde{x}} \subseteq \Pi_{\overline{X}_x^{\mr{log}}}$ (\text{resp.}, 
$D_{\tilde{y}} \subseteq \Pi_{\overline{Y}_y^{\mr{log}}}$) a specific decomposition group of the diagonal cusp $\tilde{x}^{\mr{log}}$  (\text{resp.}, $\tilde{y}^{\mr{log}}$)
[cf. the discussion following Remark 2.1.1].

 Let 
\[ \alpha:\Pi_{\overline{X}_x^{\mr{log}}} \stackrel{\sim}{\longmigi} \Pi_{\overline{Y}_y^{\mr{log}}} \]
 be a Frobenius-preserving isomorphism [cf. Definition 3.7 (i)] which maps $D_{\tilde{x}}$ onto $D_{\tilde{y}}$.
Let us denote by $\overline{\alpha} : \Pi_X \isom \Pi_Y$
the isomorphism obtained by passing to the quotients $\Pi_{\overline{X}_x^{\mr{log}}} \migisurj \Pi_X$, $\Pi_{\overline{Y}_y^{\mr{log}}} \migisurj \Pi_Y$. 
Let us denote by
$D_x \subseteq \Pi_X$ (resp., $D_y \subseteq \Pi_Y$) 
the decomposition group of $x$ (resp., the decomposition group of $y$)  determined by the image of $D_{\tilde{x}}$ in $\Pi_X$ (resp., as the image of $D_{\tilde{y}}$ in $\Pi_Y$) via the quotient $\Pi_{\overline{X}_x^{\mr{log}}} \migisurj \Pi_X$ (resp., $\Pi_{\overline{Y}_y^{\mr{log}}} \migisurj \Pi_Y$). 

Then there exists an isomorphism
\[ \alpha_2:\Pi_{X_2} \stackrel{\sim}{\longmigi} \Pi_{Y_2} \]
which is uniquely determined up to composition with an inner automorphism (of either the domain or codomain) by the condition that
it
 is compatible with the natural switching automorphisms [cf. the discussion following Remark 2.1.1] and with the specific decomposition groups associated to the respective diagonal divisors determined by $D_{\tilde{x}}$, 
$D_{\tilde{y}}$  [cf. Lemma 2.4 (ii)], 
which fits into the following commutative square
\[ \hspace{45mm} \begin{CD}
\Pi_{X_2} @> \alpha_2 >> \Pi_{Y_2}
\\
@V p_1^\Pi VV @VV p_1^\Pi V 
\\
\Pi_X @> \overline{\alpha} >> \Pi_Y,
\end{CD} \hspace{40mm} (**) \]
and induce $\alpha$ upon restriction to the inverse images (via the vertical arrows of (**)) of $D_x\subseteq \Pi_X$ and $D_y\subseteq \Pi_Y$. 
\et
%--------------------------------------------------------------[begin proof]-----------------
%\vspace{2mm}
\bpf
Let us fix specific choices of $(C_{i_1}^X)$, $(C_{p_1}^X)$, $(C_{i_1}^Y)$, $(C_{p_1}^Y)$.
By applying Proposition 3.8 to these choices and the choices of $(C_{\tilde{x}}^X)$, $(C_{\tilde{y}}^Y)$, $(C^{\overline{X}^{\mr{log}}_x \isom \overline{Y}^{\mr{log}}_y}_{x,y})$ given by hypothesis, we obtain a commutative diagram as follows:  
\bcd
\Pi_{\overline{X}_x^{\mr{log}}} @> \mr{Int}_X^\Pi \circ i_1^\Pi >> \Pi^{\mr{Lie}}_{X_2} @> s_{X}^{\Pi^{\mr{Lie}}} >> \Pi^{\mr{Lie}}_{X_2}
\\
@V \alpha VV  @V \alpha_2^{\Pi^{\mr{Lie}}} VV @V \breve{\alpha}_2^{\Pi^{\mr{Lie}}} VV
\\
\Pi_{\overline{Y}_y^{\mr{log}}} @> \mr{Int}_Y^\Pi \circ i_1^\Pi >> \Pi^{\mr{Lie}}_{Y_2} @> s_{Y}^{\Pi^{\mr{Lie}}} >> \Pi^{\mr{Lie}}_{Y_2}
\ecd
Now observe that by the various constructions involved,
$s_X^{\Pi^{\mr{Lie}}} \circ s_X^{\Pi^{\mr{Lie}}} = \mr{id}_{\Pi^{\mr{Lie}}_{X_2}
}$, and 
$s_X^{\Pi^{\mr{Lie}}} \circ \mr{Int}_X^\Pi \circ i_1^\Pi$ coincides with 
$\mr{Int}_X^\Pi \circ i_2^\Pi$ for some $i_2^\Pi :\Pi_{\overline{X}_x^{\mr{log}}} \migi \Pi^{}_{\overline{X}^{\mr{log}}_2}$ (within the conjugacy class of homomorphisms determined by 
$i_2^\Pi$)
induced by $i_2 : \overline{X}_x^{\mr{log}} \migi \overline{X}^{\mr{log}}_2$.
Thus, it follows from Proposition 1.4 (ii) that $(\mr{Int}_X^\Pi \circ i_1^\Pi )( \Pi_{\overline{X}_x^{\mr{log}}})$ and
  $(s_{X}^{\Pi^{\mr{Lie}}} \circ \mr{Int}_X^\Pi \circ i_1^\Pi)(\Delta_{X_{2/1}}^1)$ 
generate $\Pi_{X_2}$,
and that $\Pi_{X_2}$ is preserved by the action of $s_X^{\Pi^{\mr{Lie}}}$.
Similarly, $\Pi_{Y_2}$ is generated by $(\mr{Int}_Y^\Pi \circ i_1^\Pi)(\Pi_{\overline{Y}_y^{\mr{log}}})$ and $(s_{Y}^{\Pi^{\mr{Lie}}} \circ \mr{Int}_Y^\Pi \circ i_1^\Pi)(\Delta_{Y_{2/1}}^1)$,  and $\Pi_{Y_2}$ is preserved by the action of  $s_Y^{\Pi^{\mr{Lie}}}$.
Therefore, since the above diagram is
commutative, $\alpha_2^{\Pi^{\mr{Lie}}}$ maps $\Pi_{X_2}$ onto $\Pi_{Y_2}$.
Thus, the restriction $\alpha_2$ of $\alpha_2^{\Pi^{\mr{Lie}}}$ to $\Pi_{X_2}$ makes the diagram $(**)$ commute and is compatible with the switching automorphisms. 
Since the specific inertia subgroup of $\Pi_{X_2}$ associated to the diagonal divisor determined by $D_{\tilde{x}}$ is the image of $I_{\tilde{x}} \subseteq \Pi_{\overline{X}_x^{\mr{log}}}$ via $\mr{Int}_X^\Pi \circ i_1^\Pi$ (cf. Lemma 2.4 (ii)), 
the isomorphism $\alpha_2$, which is an extension of the isomorphism $\alpha$, is compatible with the corresponding specific decomposition groups associated to the respective diagonal divisors.
This completes the proof of the existence assertion.

Next, we consider uniqueness.
Let $\dot{\alpha}_2$, $\ddot{\alpha}_2$ $: \Pi_{X_2} \isom \Pi_{Y_2}$ be isomorphisms both of which make the diagram $(**)$ commute and induce $\alpha |_{\Delta_{X_{2/1}}^1}$ (i.e., the restriction of $\alpha$ to $\Delta_{X_{2/1}}^1 $) upon restriction to the kernels of the vertical arrows of $(**)$. 
Then $\dot{\alpha}_2^{-1} \circ \ddot{\alpha}_2$
determines an automorphism of the exact sequence
\[ 1 \longmigi \Delta_{X_{2/1}}^1 \stackrel{i_1^\Pi}{\longmigi} \Pi_{X_2} \stackrel{p_1^\Pi}{\longmigi} \Pi_X \longmigi 1 \]
which induces the identity automorphisms on $\Delta_{X_{2/1}}^1$ and  $\Pi_X$. This implies that $\dot{\alpha}_2^{-1} \circ \ddot{\alpha}_2$ is the identity morphism (cf. the last paragraph of ``{\bf Topological Groups}" in \S\,0). 
\epf
%------------------------------------------------------------------[end theorem]---------------
%-------------------------------------------------------------------------------[begin corollary]---------
\bco\leavevmode\\
 \ \ \ 
Let $X$ (\text{resp.}, $Y$) be a hyperbolic curve over a finite field $K$ (\text{resp.}, $L$),
$x, x'$ $K$-rational points of $\overline{X}$ (\text{resp.}, $y, y'$ $L$-rational points of $\overline{Y}$).
 Let 
\[ \alpha:\Pi_{\overline{X}_x^{\mr{log}}} \longmigi \Pi_{\overline{Y}_y^{\mr{log}}} \]
 be a Frobenius-preserving isomorphism
 such that the decomposition groups of $\tilde{x}$ and $\tilde{y}$
 (which are well-defined up to conjugacy)
 correspond via $\alpha$.
 Suppose that
the isomorphism $\overline{\alpha} : \Pi_X \isom \Pi_Y$ induced by passing to the quotients $\Pi_{\overline{X}_x^{\mr{log}}} \migisurj \Pi_X$, $\Pi_{\overline{Y}_y^{\mr{log}}} \migisurj \Pi_Y$ maps the conjugacy class of the decomposition group of $x'$ to the conjugacy class of the decomposition group of $y'$. 

 Then there exists a Frobenius-preserving isomorphism 
\[ \alpha':\Pi_{\overline{X}_{x'}^{\mr{log}}} \longmigi \Pi_{\overline{Y}_{y'}^{\mr{log}}} \]
which is uniquely determined up to composition with an inner automorphism (of either the domain or codomain) by the condition that
 it  induces $\overline{\alpha}$ upon passing to the respective quotients and maps the conjugacy class of the decomposition group of the diagonal cusp $\tilde{x'}$ to the conjugacy class of the decomposition group of the diagonal cusp $\tilde{y'}$. 
\eco
%---------------------------------------------------------------------------[begin proof]----------------
%\vspace{2mm}
\bpf 
The existence assertion follows from Theorem 3.9 and the fact that if $D_{x'} \subseteq \Pi_X$, $D_{y'} \subseteq \Pi_Y$ denote the decomposition groups of $x'$, $y'$ respectively, then we have natural isomorphisms $\Pi_{\overline{X}_{x'}^{\mr{log}}} \cong D_{x'} \times_{\Pi_X} \Pi_{X_2}$, $\Pi_{\overline{Y}_{y'}^{\mr{log}}} \cong D_{y'} \times_{\Pi_Y} \Pi_{Y_2}$. 

Next, we consider the uniqueness assertion.
Let  $\dot{\alpha}', \ddot{\alpha}' : \Pi_{\overline{X}_{x'}^{\mr{log}}} \isom \Pi_{\overline{Y}_{y'}^{\mr{log}}}$ be Frobenius-preserving isomorphisms both of which 
 induce $\overline{\alpha}$ upon passing to the respective quotients and map some {\it specific} decomposition group of the diagonal cusp $\tilde{x'}$ to the {\it same} decomposition group of the diagonal cusp $\tilde{y'}$. 
Write $\beta := (\dot{\alpha}')^{-1} \circ  \ddot{\alpha}' \in \mr{Aut}(\Pi_{\overline{X}_{x'}^{\mr{log}}})$.
Then it follows from the existence portion of Theorem 3.9
that $\beta$ induces an element  $\beta_2 \in \mr{Aut}(\Pi_{X_2})$ which induces the identity morphism of $\Pi_{X^{\times 2}}$ upon passing to the natural quotient $\Pi_{X_2} \migisurj \Pi_{X^{\times 2}}$.
Note that $\beta_2$ defines an element $[\beta_2]\in\mr{Out}^{\mr{FC}}(\Delta_{X_2})$.
Moreover,  since $\beta_2$ induces the identity morphism of $\Pi_{X^{\times 2}}$, it follows that
$[\beta_2]$ maps to the identity element of $\mr{Out}(\Delta_{X})$ (cf. ~\cite{HSHMZK}  for the definition of and results concerning to ``$\mr{Out}^{\mr{FC}}$").
But $\mr{Out}^{\mr{FC}}(\Delta_{X_2}) \migi \mr{Out}(\Delta_{X})$ is injective (cf., e.g., ~\cite{HSHMZK}, Theorem A), so we have $[\beta_2] = 1$ i.e., the restriction of $\beta_2$ to $\Delta_{X_2}$ coincides with  an inner automorphism $\mr{Inn}(b)$ determined by an element $b$ of  $\Delta_{X_2}$. 
By the construction of $\beta_2$, $\mr{Inn}(b)$ (preserves the subgroup $\Delta_{X_{2/1}}^1$ of $\Delta_{X_2}$ and) induces the identity morphism of $\Delta_{X}$ upon passing to the quotient $\Delta_{X_2} \migisurj \Delta_{X_2} / \Delta_{X_{2/1}}^1 \cong \Delta_X$.
Since $\Delta_X$ is center-free (cf. Proposition 1.4 (iii)), we thus conclude that $b$ maps to the identity element of $\Delta_X$ via $\Delta_{X_2} \migisurj  \Delta_X$.
In particular, $b$ is an element of $\Pi_{\overline{X}_{x'}^{\mr{log}}}$. 
 Thus we have  two automorphisms $\beta$, $\mr{Inn}(b)$ on $\Pi_{\overline{X}_{x'}^{\mr{log}}}$  which coincide upon passing to the quotient $\Pi_{\overline{X}_{x'}^{\mr{log}}} \migisurj D_{x'} \subseteq \Pi_X$ as well as upon the restriction to $\Delta_{X_{2/1}}^1 \subseteq \Pi_{\overline{X}_{x'}^{\mr{log}}}$. This implies that $\beta = \mr{Inn}(b)$ (cf. the last paragraph of ``{\bf Topological Groups}" in \S 0), hence  completes the proof of the uniqueness assertion. 
\epf
%\vspace{1mm}
%----------------------------------------------------------------------------[end corollary]------------
%-------------------------------------------------------------------------[begin remark]---------------
\begin{rema}
\leavevmode\\
 \ \ \ Any Frobenius-preserving isomorphism is quasi-point-theoretic (cf. ~\cite{TAM}, Corollary 2.10, Proposition 3.8; ~\cite{MZKcomb1}, Remark 10, (iii)), i.e., induces a bijection between the sets of decomposition groups of the points of $\overline{X}, \overline{Y}$.
Therefore, 
in the statement of Corollary 3.10, given a closed point $x''$ of $\overline{X}$, there always {\it exists} a closed point $y''$ of $\overline{Y}$ which corresponds, at the level of conjugacy classes of decomposition groups, to $x''$ via $\overline{\alpha}$ (but this choice is {\it not necessarily unique}!).
\end{rema}
%-------------------------------------------------------------------------------[end remark]-------------
\leavevmode\\
%%%%%%%%%%%%%%%%%%%%%%%%%%%%%%%%%%%%%%%%%%%%%%%%
%%%%%%%%%%%%%%%%%%%%%%%%%%%----[begin section]----%%%%%%%%%%
\section{Cuspidalization Problems for hyperbolic curves} \leavevmode\\
 %-----------------------------------------------------------------------------------------------------------
%------------------------------------------------------------------------------------------------
 \ \ \ In this last section, we apply Theorem 3.9 to obtain group-theoretic constructions of the cuspidalization of a hyperbolic curve at a point infinitesimally close to a cusp (cf. Theorem 4.3), as well as of arithmetic fundamental groups of configuration spaces of arbitrary dimension (cf. Theorem 4.4). 
\\[3mm]  
 \ \ \ We maintain the notation and set-up of the discussion at the beginning of Section 3.
 Moreover, until the end of Theorem 4.3,
 we shall assume that $X$ is affine (i.e., $r>0$),
 and that $x$ is a split cusp of $X$, i.e., $x \in \overline{X}(K) \setminus X(K)$.
As discussed following Remark 2.1.1, the major and minor cuspidal components $\overline{X}^{\mr{log}'}$, $\overline{\mbP}_X^{\mr{log}'}$ at $x$, together with the nexus $\nu_x^{\mr{log}}$ at $x$, determine {\it strict}
 (cf. ~\cite{ILL}, 1.2) closed sub-log schemes of $\overline{X}_x^{\mr{log}}$.  These closed sub-log schemes determine subgroups well-defined up to conjugacy
\[ \Pi_{\overline{X}^{\mr{log}'}}, \ \Pi_{\overline{\mbP}_X^{\mr{log}'}}, \ \Pi_{\nu_x^{\mr{log}}} \subseteq \Pi_{\overline{X}_x^{\mr{log}}} \]  
--- which we shall refer to, respectively, as the \textit{major verticial}, \textit{minor verticial}, and \textit{nexus subgroups} (cf. [14], Definition 1.4) ---. 

%--------------------------------------------------------------------------[begin lemma]------------------
%\vspace{3mm}
\ble \leavevmode\\
 \ \ \ 
 Write
  \[ D_x := \mr{Im}(\Pi_{\overline{X}^{\mr{log}}_x} \stackrel{p^\Pi_1 \circ i^\Pi_1}{\migi} \Pi_X).\]
 (Thus, $D_x \subseteq \Pi_X$ is a specific decomposition group of $x$, i.e., well-defined without any conjugacy indeterminacies.)  Then:
 
 \begin{itemize}
 \item[(i)] For any choice of a specific major verticial subgroup $\Pi_{\overline{X}^{\mr{log}'}} \subseteq \Pi_{\overline{X}_x^{\mr{log}}}$ ,
the composite morphism
\[ \Pi_{\overline{X}^{\mr{log}'}} \stackrel{}{\longmigi} \Pi_{\overline{X}^{\mr{log}}_x} \stackrel{(p^\Pi_1 \circ i^\Pi_1, p^\Pi_2 \circ i^\Pi_1)}{\longmigi} D_x \times_{G_K} \Pi_X  \]
 is an isomorphism.
(In particular,
the major verticial subgroups may be thought of
as defining sections of the natural surjection $\Pi_{\overline{X}_x^{\mr{log}}} \migisurj \Pi_X \times_{G_K} D_x$.)
Moreover, the inverse of this isomorphism maps the subgroup $D_x \times_{G_K} D_x \subseteq D_x \times_{G_K} \Pi_X$ to the nexus subgroup $\Pi_{\nu_x^{\mr{log}}} \subseteq \Pi_{\overline{X}^{\mr{log}'}}$.

\item[(ii)] In a similar vein, 
 let $\overline{\mbP}_K^{\mr{log}}$ be the 1-st log configuration space associated to a tripod $\mbP_K$ over $K$ (cf. Definition 1.1 (ii)).
 Then  for any choice of a specific minor verticial subgroup $\Pi_{\overline{\mbP}_X^{\mr{log}'}} \subseteq \Pi_{\overline{X}_x^{\mr{log}}}$,
the composite morphism
\[ \Pi_{\overline{\mbP}_X^{\mr{log}'}} \stackrel{}{\longmigi} \Pi_{\overline{X}^{\mr{log}}_x} \stackrel{(p^\Pi_\mbP,  p^\Pi_1 \circ i^\Pi_1)}{\longmigi} \Pi_{\overline{\mbP}^{\mr{log}}_K}  \times_{G_K} D_x   \]
--- where $p^\Pi_\mbP$ denotes
the homomorphism $\Pi_{\overline{X}^{\mr{log}}_x} \migi \Pi_{\overline{\mbP}^{\mr{log}}_K}$
(well-defined up to conjugation) induced by the natural morphism $\overline{X}^{\mr{log}}_x \migi \overline{\mbP}^{\mr{log}}_K$ given by contracting  $\overline{X}$ ($\subseteq \overline{X}_x$) 
to $\nu_x$ --- 
 is an isomorphism.
\end{itemize}
\ele
%-------------------------------------------------------------------------[begin proof]-----------------
%\vspace{2mm}
\bpf
We shall only consider assertion (i) since assertion (ii) follows from a similar argument.
Let us consider the commutative diagram of natural morphisms of log schemes
\bcd
\nu_x^{\mr{log}} @>>> x^{\mr{log}} \times_K x^{\mr{log}}
\\
@VVV @VVV
\\
\overline{X}^{\mr{log}'}  @>>> \overline{X}^{\mr{log}} \times_K x^{\mr{log}}
\ecd
--- where the horizontal arrows are the strict closed immersions ---.
Now recall that: (a) k\'{e}t coverings may be constructed by means of descent with respect to (non-logarithmic!)  \'{e}tale morphisms; (b) restriction from a henselian trait to its closed point induces an equivalence between the respective categories of k\'{e}t coverings (cf. ~\cite{ILL}). 
Since the bottom horizontal arrow $\overline{X}^{\mr{log}'} \migi \overline{X}^{\mr{log}} \times_K x^{\mr{log}}$ in the above diagram is an isomorphism on the respective complements of the images of the horizontal arrows in the above diagram,
it suffices (by (a), (b)) to verify that
 the induced morphism between the log inertia groups of $\nu_x^{\mr{log}}$ and $x^{\mr{log}} \times_K x^{\mr{log}}$ 
(i.e., $\mr{Ker}(\Pi_{\nu_x^{\mr{log}}} \migisurj G_K)$ and $\mr{Ker}(\Pi_{x^{\mr{log}} \times_K x^{\mr{log}}} \migisurj G_K)$)
is an isomorphism
 (cf. ~\cite{ILL}, 4.7 for the terminology ``log inertia subgroup"). 
Fix a chart, modeled on $\mbN$, of $x^{\mr{log}}$ (i.e., roots of a local uniformizer at $x$ in $\overline{X}$).
Then such a chart determines charts, modeled on $\mbN \oplus \mbN$,
of
$x^{\mr{log}} \times_K x^{\mr{log}}$ 
and
$\nu_x^{\mr{log}}$.
By using these charts, one verifies easily that 
the homomorphism of monoids induced by the morphism $\nu_x^{\mr{log}} \migi x^{\mr{log}} \times_K x^{\mr{log}}$ may be expressed as follows:
\[ \mbN \,\oplus \,\mbN \longmigi \mbN \,\oplus \,\mbN \]
\[ \ \ \ \ ( a, b ) \ \mapsto \ ( a+b , b ). \]
Then, by applying the functor $\mr{Hom}(( \ \_ \ )^{\mr{gp}}, \mbZ_l (1))$ to this morphism of monoids, one verifies immediately that the induced morphism of log inertia groups between $\nu_x^{\mr{log}}$ and $x^{\mr{log}} \times_K x^{\mr{log}}$ is an isomorphism.
\epf
%---------------------------------------------------------------------------[end lemma]------------------
%----------------------------------------------------------------------------[begin lemma]----------------
%\vspace{3mm}
\ble \leavevmode\\
 \ \ \ Suppose that we fix a choice of a nexus subgroup $\Pi_{\nu_x^{\mr{log}}} \subseteq \Pi_{\overline{X}_x^{\mr{log}}}$ among its various $\Pi_{\overline{X}_x^{\mr{log}}}$-conjugates.
Then:
\be
\arrangeVa\item[(i)]
There exists a unique pair of inclusions
\[ \Pi_{\overline{X}^{\mr{log}'}} \subseteq \Pi_{\overline{X}_x^{\mr{log}}}, \ \ \ \Pi_{\overline{\mbP}_X^{\mr{log}'}} \subseteq \Pi_{\overline{X}_x^{\mr{log}}} \]
(among their various $\Pi_{\overline{X}_x^{\mr{log}}}$-conjugates) both of which contain $\Pi_{\nu_x^{\mr{log}}} \subseteq \Pi_{\overline{X}_x^{\mr{log}}}$ . 
\arrangeVa\item[(ii)]
 The  inclusions $\Pi_{\nu_x^{\mr{log}}} \subseteq \Pi_{\overline{X}^{\mr{log}'}} \subseteq \Pi_{\overline{X}_x^{\mr{log}}}$, $\Pi_{\nu_x^{\mr{log}}} \subseteq \Pi_{\overline{\mbP}_X^{\mr{log}'}} \subseteq \Pi_{\overline{X}_x^{\mr{log}}}$ obtained in (i) make  the diagram
\bcd
\Pi_{\nu_x^{\mr{log}}} @>>> \Pi_{\overline{\mbP}_X^{\mr{log}'}}
\\
@VVV @VVV
\\
\Pi_{\overline{X}^{\mr{log}'}} @>>> \Pi_{\overline{X}_x^{\mr{log}}}
\ecd
commute and  co-cartesian in the category of profinite groups equipped with an augmentation to  $G_{K}$ whose kernel is pro-$l$.
\ee
\ele
%-------------------------------------------------------------------------------[begin proof]------------
%\vspace{2mm}
\bpf
Assertion (i) (respectively, (ii)) follows immediately from ~\cite{MZKcomb2},
Proposition 1.5, (ii) (respectively 1.5, (iii)).
\epf
%------------------------------------------------------------------------------[end lemma]---------------
%------------------------------------------------------------------------------------------------
Next, we turn to the proof of Theorem B. 
Theorem 4.3 given below may be regarded as a slightly weakened version of Theorem B (as stated in the Introduction).  This weakened version, however, will be sufficient to prove Theorem 4.4 below (which corresponds {\it precisely} to Theorem C in the Introduction).  Moreover,  one may conclude Theorem B (as stated in the Introduction) from Theorem C (cf.  Remark 4.4.1).
On the other hand, if we did {\it not} restrict our attention, in the statement of Theorem 4.3, to this slightly weakened version of Theorem B, then it would have been necessary to (essentially) {\it repeat}, in our proof of Theorem 4.4 below, arguments already applied in the proof of Theorem 4.3.

%-------------------------------------------------------------------------------[begin Theorem]-------
%\vspace{3mm}
\bt \leavevmode\\
 \ \ \ 
 Let $X$ (resp., $Y$) be an affine hyperbolic curve over a finite field $K$ (resp., $L$),
$x$ (resp., $y$) a $K$-(resp., $L$-)rational point of $\overline{X} \setminus X$
(resp., $\overline{Y} \setminus Y$).
Let 
\[ \alpha : \Pi_{X } \stackrel{\sim}{\longmigi} \Pi_{Y } \]
 be a Frobenius-preserving isomorphism such that the decomposition groups of $x$ and $y$ (which are well-defined up to conjugacy) correspond via $\alpha$.
In the following, we shall apply the notational conventions introduced in the discussion following Remark 2.1.1.

Then there exist finite extensions $\dot{K}$ of $K$ and $\dot{L}$ of $L$ and
an isomorphism 
\[ \dot{\alpha}_{x,y} : \Pi_{\overline{X}_x^{\mr{log}}}\times_{G_K} G_{\dot{K}} \stackrel{\sim}{\longmigi} \Pi_{\overline{Y}_y^{\mr{log}}} \times_{G_L} G_{\dot{L}}  \]
which is uniquely determined up to composition with an inner automorphism
(of either the domain or codomain)
by the condition that it maps the 
conjugacy class of the
decomposition group
of $\tilde{x}$
to 
the  conjugacy class of the decomposition group
 of $\tilde{y}$
and induces $\alpha |_{ \Pi_{X \times_K \dot{K}}} :  \Pi_{X \times_K \dot{K}} \isom \Pi_{Y \times_L \dot{L}}$ upon  passing to the quotients $\Pi_{\overline{X}_x^{\mr{log}}} \times_{G_K} G_{\dot{K}}  \migisurj \Pi_{X \times_K \dot{K}}$, $\Pi_{\overline{Y}_y^{\mr{log}}} \times_{G_L} G_{\dot{L}} \migisurj \Pi_{Y \times_L \dot{L}}$.

\et
%------------------------------------------------------------------------------[begin proof]----------
%\vspace{2mm}
\bpf
The asserted uniqueness follows immediately from the uniqueness portion of Corollary 3.10.
 Next, we shall consider the existence assertion.
First, observe that there exists a connected finite \'etale covering $f:\dot{Z}\migi X$, where
$\dot{Z}$ is a hyperbolic curve over a finite extension field $\dot{K}$ of $K$ whose  (smooth) compactification admits at least two distinct $\dot{K}$-rational points $z$, $z'$  lying over $x$ at which $f$ is {\it unramified}.
Indeed, this follows immediately from the well-known structure of $\Delta_X$.
In the following, we  shall, for simplicity,  replace $\dot{K}$ by $K$  (i.e., assume that the base fields of $X$ and $\dot{Z}$ coincide).

Write
 $Z$ for the partial (smooth) compactification of $\dot{Z}$ at $z'$ 
 (i.e., 
 a unique open subscheme $Z$ of the smooth compactification of $\dot{Z}$ containing $\dot{Z}$ and satisfying that $Z \setminus \dot{Z} = \{ z'\}$)
% the hyperbolic curve $Z$ obtained from $\dot{Z}$ by adding the point $z'$ in the smooth compactification of $\dot{Z}$)
 and $\overline{Z}_{z}^{\mr{log}}$ for the cuspidalization of $Z$ at $z$.
 Thus, the underlying scheme $\overline{Z}_z$ of $\overline{Z}_{z}^{\mr{log}}$ is {\it proper}.
Denote by
\[ \overline{Z}^{\mr{log}'}, \ \ \ 
\overline{\mbP}_Z^{\mr{log}'}, \ \ \ 
\nu_{z}^{\mr{log}}\]
the major and minor cuspidal components and the nexus of $\overline{Z}_{z}^{\mr{log}}$ at $z$, respectively (cf. the discussion at the beginning of the present \S 4).    
Let us fix  specific choices of the decomposition groups $\dot{D}_{z} \subseteq \Pi_{\dot{Z}}$ of $z$ and $D_{x} \subseteq \Pi_{X}$ of $x$ such that $D_x \cap \Pi_{\dot{Z}}  =  \dot{D}_z$.
Denote by $D_z$ the image of $\dot{D}_z$ via the quotient $\Pi_{\dot{Z}} \migisurj \Pi_Z$ (which may be considered as the decomposition group of $z$ in $\Pi_Z$).
Thus, the natural inclusion $\dot{D}_z\subseteq D_x$ is in fact an equality $\dot{D}_z=D_x$, and we have a natural isomorphism $\dot{D}_z\isom D_z$.
By applying Corollary 3.10 (cf. also Theorem 3.9) to the hyperbolic curve $\dot{Z} = Z \setminus \{ z' \}$ together with the $K$-rational points $z$ and $z'$,
we may reconstruct, group-theoretically from $\Pi_{\dot{Z}}$, the profinite group $\Pi_{\overline{Z}_{z}^{\mr{log}}}$ together with its natural augmentation to $D_z$.
Also, by ~\cite{MZKcomb1}, Corollary 2.7 (iii),  we may reconstruct, group-theoretically from the natural augmentation $\Pi_{\overline{Z}_{z}^{\mr{log}}} \migisurj D_z$,  
 the conjugacy classes of the major verticial, minor verticial and nexus subgroups of $\Pi_{\overline{Z}_{z}^{\mr{log}}}$ associated to the cuspidalization at $z$.
Now let us 
fix  {\it specific choices} of  the major verticial, minor verticial and nexus subgroups of $\Pi_{\overline{Z}_{z}^{\mr{log}}}$
\[  \Pi_{\overline{Z}^{\mr{log}'}}, \ \ \  \Pi_{\overline{\mbP}_Z^{\mr{log}'}}, \ \ \ \Pi_{\nu_{z}^{\mr{log}}} \]
such that:
(a)
the subgroup $\Pi_{\nu_{z}^{\mr{log}}} \subseteq \Pi_{\overline{Z}_{z}^{\mr{log}}}$ maps, via the natural morphism $\Pi_{\overline{Z}_{z}^{\mr{log}}} \migisurj  \Pi_Z$,
onto the subgroup $ D_z$;
(b)
$\Pi_{\nu_{z}^{\mr{log}}} \subseteq \Pi_{\overline{Z}^{\mr{log}'}} \subseteq \Pi_{\overline{Z}_{z}^{\mr{log}}}$ and $\Pi_{\nu_{z}^{\mr{log}}} \subseteq \Pi_{\overline{\mbP}_Z^{\mr{log}'}} \subseteq  \Pi_{\overline{Z}_{z}^{\mr{log}}}$.
(These choices are possible by virtue of Lemmas 4.1(i), 4.2.)
If we denote by $\overline{\mbP}^{\mr{log}}_K$ the 1-st log configuration space associated to a tripod $\mbP_K$ over $K$, then we obtain (cf. Lemma 4.1 (ii))
 a  composite
\[  \Pi_{\nu_{z}^{\mr{log}}}  \stackrel{}{\longmigi} \Pi_{\overline{\mbP}^{\mr{log}'}_Z} \stackrel{\sim}{\longmigi} \Pi_{\overline{\mbP}^{\mr{log}}_K} \times_{G_K} D_x. \]
Here, we may regard $\Pi_{\overline{\mbP}^{\mr{log}}_K} $ as an object group-theoretically reconstructed from $\Pi_{\overline{\mbP}^{\mr{log}'}_Z}$ by thinking
of $\Pi_{\overline{\mbP}^{\mr{log}}_K}  $ as the quotient of the kernel of the natural composite augmentation $\Pi_{\overline{\mbP}^{\mr{log}'}_Z} \migisurj  D_z \migisurj G_K$ (i.e., which is naturally isomorphic to $\Delta_{\mbP_K} \times \mbZ_l(1)$) by its {\it center} (i.e., $\mbZ_l(1) $ --- cf. Proposition 1.4 (iii)). 
Also, we obtain (cf. Lemma 4.1 (i)) a  diagram of natural morphisms
\[ D_x \times_{G_K} \Pi_X \stackrel{}{\longhidari} \dot{D}_{z} \times_{G_K} \dot{D}_{z} \stackrel{\sim}{\longmigi} D_{z} \times_{G_K} D_{z} \stackrel{\sim}{\longhidari} \Pi_{\nu_{z}^{\mr{log}}}  \]
induced, by  restriction, from a diagram of natural morphisms
\[ D_x \times_{G_K} \Pi_X \ \ \hidariincl\ \ \dot{D}_{z} \times_{G_K} \Pi_{\dot{Z}} \  \ \migisurj\ \ D_{z} \times_{G_K} \Pi_{Z} \  \ \stackrel{\sim}{\hidari}\ \ \Pi_{\overline{Z}^{\mr{log}'}} .
  \]
Thus, for  suitable choices of  the subgroups   $\Pi_{\overline{X}^{\mr{log}'}}$, $\Pi_{\overline{\mbP}_X^{\mr{log}'}}$, $\Pi_{\nu_x^{\mr{log}}}$ $\subseteq \Pi_{\overline{X}_x^{\mr{log}}}$ (cf. Lemma 4.2 (i)), we obtain a natural commutative diagram:
\[ \begin{CD}
\Pi_{\overline{X}^{\mr{log}'}} @<<< \Pi_{\nu_x^{\mr{log}}} @>>> \Pi_{\overline{\mbP}_X^{\mr{log}'}}
 \\
 @V \wr V  V @V \wr V  V @V\wr V  V 
\\
D_x \times_{G_K} \Pi_X @<  << \Pi_{\nu_{z}^{\mr{log}}} @>  >> \Pi_{\overline{\mbP}^{\mr{log}}_K} \times_{G_K} D_{x},
 \end{CD} \]
where the vertical arrows are all isomorphisms by Lemma 4.1(i), (ii).
In particular,  it follows from Lemma 4.2 (ii) that $\Pi_{\overline{X}_x^{\mr{log}}}$
may be identified with the colimit of the lower horizontal sequence --- which, by the above discussion, may be {\it reconstructed group-theoretically} from the data $(\Pi_X,D_x\subseteq \Pi_X) $! --- in the above diagram.
Therefore, by comparing this diagram to the corresponding diagram for $Y$, the proof is completed. 
\epf
%--------------------------------------------------------------------------[end Theorem]------------
Next, we consider Theorem C, i.e., the cuspidalization problem for geometrically pro-$l$ fundamental groups of configuration spaces of (not necessarily proper) hyperbolic curves over finite fields.

%--------------------------------------------------------------------------[begin Theorem]-----------
%\vspace{3mm}
\bt {(cf. ~\cite{MZKcusp}, Theorem 3.1; ~\cite{HSHcusp}, Theorem 4.1)} \leavevmode\\
  \ \ \ Let $X$ (resp., $Y$) be a hyperbolic curve over a finite field $K$ (resp., $L$).
Let
\[ \alpha_1 : \Pi_X \stackrel{\sim}{\longmigi} \Pi_Y \]
be a Frobenius-preserving isomorphism. Then  for any $n \in \mbZ_{\geq 0}$, there exists an isomorphism
\[ \alpha_n : \Pi_{X_n} \stackrel{\sim}{\longmigi} \Pi_{Y_n} \]
which is uniquely determined up to composition with an inner automorphism 
(of either the domain or codomain)
by
the condition that it is
compatible with the natural respective outer actions of the symmetric group on $n$ letters and make the diagram
\bcd
\Pi_{X_{n+1}} @> \alpha_{n+1} >> \Pi_{Y_{n+1}} 
\\
@V p_i VV @VV p_i V
\\
\Pi_{X_n} @> \alpha_n >> \Pi_{Y_n}
\ecd
($i = 1, \cdots , n+1$) commute.\et
%---------------------------------------------------------------------------[begin proof]----------------
\bpf
First, we recall that the case where $n=2$ and $X$ is proper follows from ~\cite{MZKcusp}, Theorem 3.1.
Next, we consider the case where $n=2$ and $X$ is affine.
As we noted in Definition 3.7 (i),
$\alpha_1$ induces an isomorphism
\[  \alpha_0 : G_K \isom G_L \]
of profinite groups.
Now, by combining Theorem 3.9 and Theorem 4.3
together with the fact that $\alpha_1$ is quasi-point-theoretic (cf. Remark 3.10.1), 
 we conclude that $\alpha_1$ induces an isomorphism
\[ \dot{\alpha}_2 : \Pi_{X_2} \times_{G_K} G_{\dot{K}} \isom \Pi_{Y_2} \times_{G_L} G_{\dot{L}}
 \]
--- where $G_{\dot{K}} \subseteq G_K$, $G_{\dot{L}} \subseteq G_L$ denote open subgroups corresponding to certain finite extensions $\dot{K}$ of $K$ and $\dot{L}$ of $L$, respectively ---. 
If we denote by $\alpha^\Delta_2$ the restriction of $\dot{\alpha}_2$ to $ \Delta_{X_2}$,
then (cf. Theorem 3.9) $\alpha^\Delta_2$ maps onto $\Delta_{Y_2}$, i.e., determines an isomorphism 
\[  \alpha^\Delta_2 : \Delta_{X_2} \isom \Delta_{Y_2}.  \]
Let
 \[ \gamma_X : G_K \migi \mr{Out}^{\mr{FC}}(\Delta_{X_2}) \ \ \  (\text{resp.},  \gamma_Y : G_L \migi \mr{Out}^{\mr{FC}}(\Delta_{Y_2})) \]
(cf. ~\cite{HSHMZK}  for the definition and results concerning to ``$\mr{Out}^{\mr{FC}}$")
be  the morphism obtained by lifting elements of $G_K$ (resp., $G_L$) via the surjection $\Pi_{X_2} \migisurj G_K$ (resp., $\Pi_{Y_2} \migisurj G_L$)
  and considering the action of these elements by  conjugation.
Then 
 $\alpha_2^\Delta$, $\alpha_0$ give rise to two composites
 $\gamma_Y \circ \alpha_0$ and $[\alpha_2^\Delta ] \circ \gamma_X$
\[  \gamma_Y \circ \alpha_0, \ [\alpha_2^\Delta ] \circ \gamma_X : G_K \longmigi \mr{Out}^{\mr{FC}}(\Delta_{Y_2}) \]
--- where $[ \alpha_2^\Delta ]$ denotes the isomorphism $\mr{Out}^{\mr{FC}}(\Delta_{X_2}) \isom \mr{Out}^{\mr{FC}}(\Delta_{Y_2})$ that sends an element $g \in \mr{Aut}(\Delta_{X_2})$ to $\alpha_2^\Delta \circ g \circ ( \alpha_2^\Delta )^{-1} \in \mr{Aut}(\Delta_{Y_2})$ ---.
It follows from the constructions of $\alpha_0$, $\alpha_2^\Delta$
that $\gamma_Y \circ \alpha_0$ and $\ [\alpha_2^\Delta ] \circ \gamma_X$ coincide
after composing with the natural morphism
$\mr{Out}^{\mr{FC}}(\Delta_{Y_2}) \migi \mr{Out}(\Delta_{Y})$.
On the other hand, since
$\mr{Out}^{\mr{FC}}(\Delta_{Y_2}) \migi \mr{Out}(\Delta_{Y})$ is injective
(cf., e.g., [5], Theorem A),
we conclude that $\gamma_Y \circ \alpha_0 = [\alpha_2^\Delta ] \circ \gamma_X$.
Therefore,
by applying the natural isomorphisms $\Pi_{X_2} \cong \Delta_{X_2} \stackrel{\mr{out}}{\rtimes} G_K$ and
$\Pi_{Y_2} \cong \Delta_{Y_2} \stackrel{\mr{out}}{\rtimes} G_L$,
we obtain an isomorphism 
$\Pi_{X_2} \cong \Pi_{Y_2}$, which satisfies the required uniqueness and compatibility properties (cf. the construction of $\dot{\alpha}_2$; Theorem 3.9). This completes the proof of the assertion in the case where $n=2$ and $X$ is affine.

Finally, the assertion in the case $n \geq 3$ follows from an inductive argument on $n$ applied to an argument similar to the argument of the above discussion. 
Indeed, consider the natural exact sequence
\[ 1 \longmigi \Delta_{(X \times_K \overline{K} \setminus \{ x \} )_{n-1}} \longmigi \Pi_{X_n} \stackrel{q_j^\Pi}{\longmigi} \Pi_X \longmigi 1\]
(which induces an isomorphism $\Pi_{X_n} \cong \Delta_{(X \times_K \overline{K} \setminus \{ x \} )_{n-1}} \stackrel{\mr{out}}{\rtimes} \Pi_X$), where $x$ denotes a $\overline{K}$-rational point of $X $, and $q_j^\Pi$ denotes the morphism induced by the projection $X_n \migi X$ to the $j$-th factor.
Since the natural morphism $\mr{Out^{FC}}(\Delta_{(X \times_K \overline{K} \setminus \{ x \} )_{n-1}}) \migi \mr{Out^{FC}}(\Delta_{(X \times_K \overline{K} \setminus \{ x \} )_{n-2}})$ is {\it injective} (cf. [5]),
 we may carry out a similar argument to the above discussion by replacing 
$G_K$ by $\Pi_{X}$ and $\Delta_{X_2}$ by $\Delta_{(X \times_K \overline{K} \setminus \{ x \} )_{n-1}}$.
Hence, for $j= 1, \cdots , n$,  we obtain an isomorphism 
$\alpha^j_{n} :  \Pi_{X_{n}} \isom \Pi_{Y_{n}} $ 
that fits into a
commutative diagram
\bcd
\Pi_{X_{n}} @> \alpha^j_{n} >> \Pi_{Y_{n}} 
\\
@V p_i^\Pi VV @VV p_i^\Pi V
\\
\Pi_{X_{n-1}} @> \alpha_{n-1} >> \Pi_{Y_{n-1}}.
\ecd
for  $i= 1, \cdots , n-1$.
But it follows from the induction hypothesis (concerning to the asserted uniqueness), together with the {\it injectivity} applied above, that  the 
$\alpha^j_{n}$'s coincide, for $j=1,\cdots,n$, up to composition with an inner automorphism, and that the asserted uniqueness and compatibility with symmetric group actions for $n$ are satisfied.
\epf
%-------------------------------------------------------------------------------[end theorem]---------

%---------------------------------------------------------------------------------[begin remark]----------
\begin{rema} \leavevmode\\
 \ \ \ As explained in the discussion preceding Theorem 4.3, one may conclude Theorem B (as stated in Introduction) directly from Theorem 4.4 as follows.
 Let $X$, $Y$, $x$, $y$ and $\alpha$ be as in the statement of Theorem B.
Then, by applying Theorem 4.4 in the case $n = 2$, we obtain a unique isomorphism
\[  \alpha_2 : \Pi_{X_2} \isom \Pi_{Y_2}  \]
 well-defined up to composition with an inner automorphism of either the domain or codomain, 
which 
fits into  two commutative diagrams as follows:
\[ \begin{CD}
\Pi_{X_2} @> \alpha_2 >> \Pi_{Y_2}
\\
@V p_1^\Pi VV @VV p_1^\Pi V 
\\
\Pi_X @> \alpha >> \Pi_Y,
\end{CD} \hspace{20mm}
\begin{CD}
\Pi_{X_2} @> \alpha_2 >> \Pi_{Y_2}
\\
@V p_2^\Pi VV @VV p_2^\Pi V 
\\
\Pi_X @> \alpha >> \Pi_Y.
\end{CD}
\]
On the other hand, we may have natural identifications
$\Pi_{\overline{X}^{\mr{log}}_x} \isom (p_2^\Pi)^{-1} (D_{x})$, $\Pi_{\overline{Y}^{\mr{log}}_y} \isom (p_2^\Pi)^{-1} (D_{y})$.
Hence the right-hand diagram above induces (since $\alpha(D_x) =D_y$) an isomorphism
\[ \alpha_{x,y}  : \Pi_{\overline{X}^{\mr{log}}_x} \isom \Pi_{\overline{Y}^{\mr{log}}_y}\]
by restricting $\alpha_2$ to the inverse images (via the vertical arrows) of $D_x\subseteq \Pi_X$ and $D_y\subseteq \Pi_Y$. 
On the other hand, it follows from [13], Corollary 2.7 (i)
that $\alpha_{x,y}$ maps the conjugacy class of the decomposition group of $\tilde{x}$ to the conjugacy class of the decomposition group of $\tilde{y}$.
Thus, the left-hand commutative diagram above induces,
by restricting the upper horizontal arrow of the diagram to the domain and codomain of $\alpha_{x,y}$,
 a commutative diagram
\[ \hspace{45mm} \begin{CD}
\Pi_{\overline{X}^{\mr{log}}_x} @> \alpha_{x,y} >> \Pi_{\overline{Y}^{\mr{log}}_y}
\\
@V p_1^\Pi VV @VV p_1^\Pi V 
\\
\Pi_X @> \overline{\alpha} >> \Pi_Y,
\end{CD} \hspace{40mm}  \]
which completes the proof of Theorem B. (The proof of uniqueness is similar to the proof of the asserted uniqueness  in Corollary 3.10.)
 \end{rema}
%-------------------------------------------------------------------------------[end remark]------------

%\vspace{3mm}
Finally, we shall conclude the paper with the following corollary.
%-------------------------------------------------------------------------------[begin corollary]-------
%\vspace{3mm}
\bco {(cf. ~\cite{HSHcusp}, Corollary 4.1) } \leavevmode\\
 \ \ \ Let $X$ (resp., $Y$) be a hyperbolic curve over a finite field $K$ (resp., $L$), and $n \in \mbZ_{\geq 0}$.
Let
\[ \alpha : \Pi_X \stackrel{\sim}{\longmigi} \Pi_Y \]
be a Frobenius-preserving isomorphism, and $x_\bullet := \{ x_1, \cdots , x_n \}$ an ordered set of distinct $K$-rational points of $X$.
 Then there exist an ordered set $y_\bullet :=\{ y_1, \cdots , y_n \}$ of distinct $L$-rational points of $Y$ and an  isomorphism
\[ \tilde{\alpha} : \Pi_{X \setminus \{ x_1, \cdots , x_n \} } \stackrel{\sim}{\longmigi} \Pi_{Y \setminus \{ y_1, \cdots , y_n \} } \]
which
is uniquely determined up to composition with an inner automorphism (of either the domain or codomain) by the condition that
it induces $\alpha$ upon passing to quotients $\Pi_{X \setminus \{ x_1, \cdots , x_n \} } \migisurj \Pi_X$, $\Pi_{Y \setminus \{ y_1, \cdots , y_n \} } \migisurj \Pi_Y$
and maps the conjugacy classes of the decomposition groups of the points in $x_\bullet$ to the conjugacy classes of the decomposition groups of the points in $y_\bullet$ in the order of numbering.
\eco
%---------------------------------------------------------------------------[begin proof]----------------
\bpf
The  existence assertion follows, by induction on $n$, from Theorem 4.4 together with the fact that any  Frobenius-preserving isomorphism between hyperbolic curves over finite fields preserves the set of decomposition groups of closed points (cf. Remark 3.10.1).
The asserted uniqueness  follows from the uniqueness asserted in Corollary 3.10, applied successively to the cuspidalizations at corresponding points of $x_\bullet$, $y_\bullet$.
\epf
%--------------------------------------------------------------------------------[end corollary]--------

%%%%%%%%%%%%%%%%%%%%%%%%%%%%%%%%%%%%%%%%%%%%%%%%%%%%%%%%%%%%%%%%%%%%%%%%%%%%%%%%%%%%%%%%%%%%%%
%---------------------------------------------------------------------------------------------

\vspace{3mm}
\begin{thebibliography}{99}
\bibitem{SGA1}
A. Grothendieck et al.,
\textit{Rev\^{e}tements \'{e}tale et groupe fondamental,}
Lecture Notes in Math., {\bf 224}, Springer, Berlin, 1971.
\bibitem{HSHlog}
Y. Hoshi, 
The exactness of the log homotopy sequence,
\textit{Hiroshima Math. J}. {\bf 39} (2009), p. 61-121.
\bibitem{HSHconf}
Y. Hoshi,
On the fundamental groups of log configuration schemes,
\textit{Math. J. Okayama Univ}. {\bf 51} (2009), p.1-26.
\bibitem{HSHcusp}
Y. Hoshi,
Absolute anabelian cuspidalizations of configuration spaces of proper hyperbolic curves over finite fields,
\textit{Publ. RIMS, Kyoto Univ}. {\bf 45} (2009), p. 661-744.
\bibitem{HSHMZK}
Y. Hoshi and S. Mochizuki,
On the combinatorial anabelian geometry of nodally nondegenarate outer representations,
\textit{Hiroshima Math. J.} {\bf 41} (2011), p. 275-342.
\bibitem{ILL}
L. Illusie,
An Overview of the Work of K. Fujiwara, K. Kato and C. Nakamura on Logarithmic Etale Cohomology,
in Cohomologies p-adiques et applications arithmetiques (2), 
\textit{Asterisque} {\bf 279} (2002), p. 271-322.
\bibitem{KATO}
K. Kato,
Logarithmic structures of Fontaine-Illusie,
in \textit{Algebraic analysis, geometry, and number theory (Baltimore, MD, 1988)}, p. 191-224,
Johns Hopkins Univ. Press, Baltimore, MD, 1989.
\bibitem{KNUD}
F. F. Knudsen, 
The projectivity of the moduli space of stable curves. II. The stacks $M_{g,r}$,
\textit{Math. Scand}. {\bf 52} (1983), p. 161-199.
\bibitem{MZKext}
S. Mochizuki,
Extending Families of Curves over Log Regular Schemes,
\textit{J. reine angew. Math}. {\bf 511} (1999), p. 319-423.
\bibitem{MZKhyp}
S. Mochizuki,
The absolute anabelian geometry of hyperbolic curves,
in \textit{Galois theory and modular forms,} p. 77-122, Dev. Math., {\bf 2}, Kluwer Acad. Publ., Boston. MA, 2004.
\bibitem{MZKsect}
S. Mochizuki,
Galois sections in Absolute Anabelian Geometry,
\textit{Nagoya Math. J.} {\bf 179} (2005), p. 17-45. 

\bibitem{MZKcusp}
S. Mochizuki,
Absolute anabelian cuspidalizations of proper hyperbolic curves,
\textit{J. Math. Kyoto Univ}. {\bf 47} (2007), p. 451-539. 
\bibitem{MZKcomb1}
S. Mochizuki,
A combinatorial version of the Grothendieck conjecture,
\textit{Tohoku Math J}. {\bf 59} (2007), p. 455-479.
\bibitem{MZKcomb2}
S. Mochizuki,
On the combinatorial cuspidalization of hyperbolic curves,
\textit{Osaka J Math}. {\bf 47}, (2010), p. 651-715.
\bibitem{MZKTAM}
S. Mochizuki and A. Tamagawa,
The algebraic and anabelian geometry of configuration spaces,
\textit{Hokkaido Math}. J. {\bf 37} (2008), p. 75-131.
\bibitem{MUM}
D. Mumford,
\textit{Abelian varieties},Tata Institute of Fundamental Reserch Studies in Mathematics, No.{\bf 5} Published for the Tata Institute of Fundamental Reserch, Bombay; Oxford University Press, London, 1970.
\bibitem{NAK}
H. Nakamura,
Galois rigidity of the \'{e}tale fundamental groups of punctured projective lines,
\textit{J. Reine Angew. Math}. {\bf 411} (1990), p. 205-216.  
\bibitem{STAB}
H. Nakamura, N. Takao and R. Ueno, 
Some stability properties of Teichmuller modular function fields with pro-$l$ weight structures,
\textit{Math. Ann}. {\bf 302} (1995), p.197-213.
\bibitem{STX}
J. Stix,
A monodromy criterion for extending curves,
\textit{Int. Math}. Res. Not. 2005.
\bibitem{TAM}
A. Tamagawa,
The Grothendieck conjecture for affine curves, 
\textit{Compositio Math}. {\bf 109} (1997), p. 135-194.

\end{thebibliography}
\end{document}